\documentclass[11pt]{amsart}
\usepackage{amssymb, latexsym}
\theoremstyle{plain}
\newtheorem{theorem}{Theorem}
\newtheorem{corollary}{Corollary}

\newtheorem{proposition}{Proposition}

\theoremstyle{remark}

\newtheorem*{Remark 1}{Remark 1}
\newtheorem*{Remark 2}{Remark 2}
\newtheorem*{Remark 3}{Remark 3}
\newtheorem*{Remark 4}{Remark 4}

\numberwithin{equation}{section}

\begin{document}

\title[The Spectrum of  Diffusion Operators]
 {Explicit and Almost Explicit Spectral Calculations for   Diffusion
 Operators}

\author{Ross G. Pinsky}
\address{Department of Mathematics\\
Technion---Israel Institute of Technology\\
Haifa, 32000\\ Israel} \email{pinsky@math.technion.ac.il}
\urladdr{http://www.math.technion.ac.il/~pinsky/}

\date{}

\begin{abstract}
The diffusion operator
$$
H_D=-\frac12\frac d{dx}a\frac d{dx}-b\frac d{dx}=-\frac12\exp(-2B)\frac d{dx}a\exp(2B)\frac d{dx},
$$
 where
$B(x)=\int_0^x\frac ba(y)dy$, defined either on $R^+=(0,\infty)$  with the Dirichlet boundary condition at $x=0$,
or on $R$, can be realized as a self-adjoint operator
with respect to the density $\exp(2Q(x))dx$. The operator is unitarily equivalent to the Schr\"odinger-type operator
$H_S=-\frac12\frac d{dx}a\frac d{dx}+V_{b,a}$, where $V_{b,a}=\frac12(\frac{b^2}a+b')$.
We obtain  an explicit criterion for the existence of a compact resolvent and
explicit formulas up to the multiplicative constant 4 for the infimum of the spectrum and for the infimum of  the
essential spectrum for these operators.
We give some applications which show in particular how $\inf\sigma(H_D)$ scales when $a=\nu a_0$ and $b=\gamma b_0$,
where $\nu$ and $\gamma$ are parameters, and $a_0$ and $b_0$ are chosen from certain classes of functions.
We also give  applications to  self-adjoint, multi-dimensional diffusion operators.
\end{abstract}
\maketitle
\section{Introduction and Statement of Results}\label{Intro}
In this paper, we give an  explicit formula up to the multiplicative constant 4
for the bottom of the spectrum and for the bottom of the essential spectrum for diffusion operators
on the half-line
$R^+=(0,\infty)$ with the Dirichlet boundary condition at 0, and for diffusion operators on the entire line.
Assuming a little more regularity,
each such operator is unitarily equivalent to a certain Schr\"odinger-type operator, so we also obtain the same information for
these latter operators.
Recall that such an  operator  possesses a compact resolvent if and only if its essential spectrum is empty, or equivalently,
if and only if the infimum of its essential spectrum is $\infty$. Thus, we
 obtain a completely explicit criterion for the existence of a compact resolvent. A diffusion operator with  a compact resolvent
 is particularly nice because its transition (sub)-probability density $p(t,x,y)$ (with respect to the reversible measure)
 can be written in the form $p(t,x,y)=\sum_{n=0}^\infty\exp(-\lambda_nt)\phi_n(x)\phi_n(y)$, where
$\{\phi_n\}_{n=0}^\infty$ is a complete, orthonormal set of eigenfunctions and $\{\lambda_n\}_{n=0}^\infty$, satisfying
$0\le\lambda_0<\lambda_1\le\lambda_2\le\cdots$, are the corresponding eigenvalues.
We give some applications of the results, which show in particular how $\inf\sigma(H_D)$ scales when $a=\nu a_0$ and $b=\gamma b_0$,
where $\nu$ and $\gamma$ are parameters, and $a_0$ and $b_0$ are chosen from certain classes of functions.
At the end of the paper, we give  applications to self-adjoint, multi-dimensional diffusion operators of the form
$-\frac12\nabla\cdot a\nabla-a\nabla Q\cdot\nabla=-\frac12\exp(-2Q)\nabla\cdot a\exp(2Q)\nabla$ on $L^2(R^d,\exp(2Q)dx)$.
The  methods  and the statements of the results are analytic, but many of the formulas and results have probabilistic
import.

We begin with the theory on the half-line, wherein lies the crux of our method. The results for the entire line follow readily from
the results for the half-line.
Let $0<a\in C^1([0,\infty))$ and  $b\in C([0,\infty))$. Define $B(x)=\int_0^x\frac ba(y)dy$.
Consider the diffusion operator with divergence-form diffusion coefficient $a$ and drift $b$
\begin{equation*}
H_D=-\frac12\frac d{dx}a\frac d{dx}-b\frac d{dx}=-\frac12\exp(-2B)\frac d{dx}a\exp(2B)\frac d{dx}
\end{equation*}
 on $R^+$
  with the Dirichlet boundary condition at $x=0$.
One can realize $H_D$ as a non-negative, self-adjoint operator on $L^2(R^+,\exp(2B)dx)$ via the Friedrichs extension of the
closure of the
nonnegative
quadratic form
$$
Q_D(f,g)=\frac12\int_0^\infty (f'ag')\exp(2B)dx,
$$
 defined for $f,g\in C^1_0(R^+)$, the space
of continuously differentiable functions with compact support on $R^+$.

Let $U_B$ denote the unitary operator from $L^2(R^+)$ to $L^2(R^+,\exp(2B)dx)$
 defined by
\begin{equation*}
U_Bf=\exp(-B)f.
\end{equation*}
Assuming that $b\in C^1(R^+)$, define $H_S=U_B^{-1}H_DU_B$. One can check that
\begin{equation*}
H_S=-\frac12\frac d{dx}a\frac d{dx}+V_{b,a},
\end{equation*}
where
\begin{equation}\label{Vb}
V_{b,a}=\frac12(\frac{b^2}a+b').
\end{equation}
Assuming that $V_{b,a}=\frac12(\frac{b^2}a+b')$ is bounded from below,
one can realize the Schr\"odinger-type operator  $H_S$ as a self-adjoint operator  on  $L^2((0,\infty))$ via the Friedrichs extension of the
closure of the semi-bounded quadratic form
$$
Q_S(f,g)=\frac12\int_0^\infty (f'ag')dx+\int_0^\infty V_{b,a}fgdx,
$$
defined for $f,g\in C^1_0(R^+)$.
Assuming in addition that $\int^\infty a(x)dx=\infty$, one can prove that
$H_S$ is in the limit-point case at $\infty$, which means in particular
that $H_S$ on $C^2_0(R^+)$ is essentially self-adjoint. (A proof in the case $a=1$ can be found in
\cite[Appendix to X.1]{RS2}. It can easily be extended to $a$ satisfying the above condition.)
  Thus, the Friedrichs extension is in fact equal to the closure
of $H_S$ on $C^2_0(R^+)$.
Note also that $U_B$ preserves the Dirichlet boundary condition.
From the above considerations, it follows that the spectra and the essential spectra of $H_D$ and $H_S$ coincide; in particular,
$H_S$ is also non-negative.

Conversely, given $a>0$,  every potential $V\ge0$ can be obtained via some $b$ as in \eqref{Vb}, and
modulo an additive constant,  every potential $V$ that is bounded from  below can be obtained via some $b$ as in \eqref{Vb}.
Indeed, let $m_V=\inf_{x\in R^+}V(x)$ and let $m_V^-=m_V\wedge0$. Since $V-m_V^-\ge0$, it is easy to show that the Riccati equation
$\frac12b'+\frac12\frac{b^2}a=V-m_V^-$
has  solutions $b$ which exist for all $x\ge0$.

The essential spectrum $\sigma_{\text{ess}}(-\frac12\frac d{dx}a\frac d{dx}+V)$ of Schr\"odinger-type operators
$-\frac12\frac d{dx}a\frac d{dx}+V$
 has been well-studied. See \cite{RS4} for the results noted below.
 For example, if $a$ is bounded and bounded from 0, and  if the potential
$V$ satisfies $\lim_{x\to\infty}V(x)=\infty$, then the operator has a compact resolvent. Thus, the spectrum consists of an increasing sequence
of eigenvalues accumulating only at infinity; in particular, $\sigma_{\text{ess}}(-\frac12\frac d{dx}a\frac d{dx}+V)=\emptyset$.
On the other hand, if $V$ is a compact (or even relatively compact) perturbation
of $-\frac12\frac d{dx}a\frac d{dx}$, which occurs in particular if $\lim_{x\to\infty}V(x)=0$, then the essential spectrum of
$-\frac12\frac d{dx}a\frac d{dx}+V$ coincides with that of $-\frac12\frac d{dx}a\frac d{dx}$;
thus $\sigma_{\text{ess}}(-\frac12\frac d{dx}a\frac d{dx}+V)=[0,\infty)$.
More generally, for arbitrary $a>0$,
 the mini-max method \cite{RS2} affords an algorithm for  arriving at
  $\inf\sigma_{\text{ess}}(-\frac12\frac d{dx}a\frac d{dx}+V)$, although
this method is mainly of theoretic import and not a practical way of calculating.

The bottom of the spectrum  of $-\frac12\frac d{dx}a\frac d{dx}+V$ is of course given by the well-known variational formula:
$$
\inf\sigma(-\frac12\frac d{dx}a\frac d{dx}+V)=\inf\frac{\int_0^\infty (\frac12a(f')^2+Vf^2)dx}{\int_0^\infty f^2dx},
$$
where the infimum is over functions $0\neq f\in C^1_0(R^+)$.
 The bottom of the spectrum of $H_D$ is also given by a variational formula:
\begin{equation}\label{variational}
\inf\sigma(H_D)=\inf\frac{\frac12\int_0^\infty a(f')^2\exp(2B)dx}{\int_0^\infty f^2\exp(2B)dx},
\end{equation}
where the infimum is over $0\neq f\in C^1_0(R^+)$.

The following theorem gives
  explicit formulas up to the multiplicative constant 4 for the bottom of the spectrum and for
the bottom of the essential spectrum of $H_D$. By the spectral
invariance, this then extends to the Schr\"odinger-type operators $H_S$.
The formulas take on two possible forms, depending on whether
$\int_0^\infty\frac1{a(x)}\exp(-2B(x))dx$ is finite or infinite.
In Remark 2 after the theorem, it is shown how the proof for the case when the integral is finite can be reduced to the case
when the integral is infinite. Remark 4 after the theorem discusses the
  probabilistic import of the theorem and of the above integral.
\begin{theorem}\label{Th}
Let $0<a\in C^1([0,\infty))$ and  $b\in C([0,\infty))$. Define
$$
B(x)=\int_0^x\frac ba(y)dy.
$$
Consider the  self-adjoint diffusion operator
$$
H_D=-\frac12\frac d{dx}a\frac d{dx}-b\frac d{dx}=
-\frac12\exp(-2B)\frac d{dx}a\exp(2B)\frac d{dx}
$$
 on
$L^2(R^+,\exp(2B)dx)$ with the Dirichlet boundary condition at 0.

\noindent If
$b\in C^1(R^+)$,
$\frac{b^2}a+b'$ is bounded from below and $\int^\infty a(x)dx=\infty$, consider also  the self-adjoint Schr\"odinger-type operator
$$
H_S=-\frac12\frac d{dx}a\frac d{dx}+\frac12(\frac{b^2}a+b')
$$
on $L^2(R^+)$ with the Dirichlet boundary condition at 0.

\noindent  If
\begin{equation}\label{intcondinfin}
\int^\infty\frac1{a(x)}\exp(-2B(x))dx=\infty,
\end{equation}
define
\begin{equation}\label{ab}
\Omega^+(b,a)=\sup_{x>0}\left(\int_0^x\frac1{a(y)}\exp(-2B(y))dy\right)\left(\int_x^\infty\exp(2B(y))dy\right)
\end{equation}
and
\begin{equation}\label{abhat}
\hat\Omega^+(b,a)=\limsup_{x\to\infty}\left(\int_0^x\frac1{a(y)}\exp(-2B(y))dy\right)\left(\int_x^\infty\exp(2B(y))dy\right).
\end{equation}
If
\begin{equation}\label{intcondfin}
\int^\infty\frac1{a(x)}\exp(-2B(x))dx<\infty,
\end{equation}
let
\begin{equation}\label{hb}
h_{b,a}(x)=\int_x^\infty\frac1{a(y)}\exp(-2B(y))dy
\end{equation}
and define
\begin{equation}\label{ab-h}
\begin{aligned}
&\Omega^+(b,a)=\sup_{x>0}\left(\int_0^xh^{-2}_{b,a}(y)\frac1{a(y)}\exp(-2B(y))dy\right)\left(\int_x^\infty h^2_{b,a}(y)\exp(2B(y))dy\right)\\
&=\sup_{x>0}\left(h_{b,a}^{-1}(x)-h_{b,a}^{-1}(0)\right)\left(\int_x^\infty h^2_{b,a}(y)\exp(2B(y))dy\right),
\end{aligned}
\end{equation}
and
\begin{equation}\label{abhat-h}
\begin{aligned}
&\hat\Omega^+(b,a)=\limsup_{x\to\infty}\left(\int_0^xh^{-2}_{b,a}(y)\frac1{a(y)}\exp(-2B(y))dy\right)\left(\int_x^\infty h^2_{b,a}(y)\exp(2B(y))dy\right)\\
&=\limsup_{x\to\infty}\left(h_{b,a}^{-1}(x)-h_{b,a}^{-1}(0)\right)\left(\int_x^\infty h^2_{b,a}(y)\exp(2B(y))dy\right).
\end{aligned}
\end{equation}
 Then
\begin{equation}\label{spec}
\frac1{8\Omega^+(b,a)}\le\inf\sigma(H_D)=\inf\sigma(H_S)\le\frac1{2\Omega^+(b,a)}
\end{equation}
and
\begin{equation}\label{essspec}
\frac1{8\hat \Omega^+(b,a)}\le\inf\sigma_{\text{ess}}(H_D)=\inf\sigma_{\text{ess}}(H_S)\le\frac1{2\hat \Omega^+(b,a)}.
\end{equation}
In particular, $H_D$ and $H_S$ possess compact resolvents if and only if $\hat \Omega^+(b,a)=0$.

\end{theorem}

\medskip

\bf\noindent Remark 1.\rm\ There does not exist a $C$ for which
$\inf\sigma(H_D)=\frac C{\Omega^+(b,a)}$, for all  drifts $b$ and all diffusion coefficients $a$.
Indeed, on the one hand, consider the case that $b(x)=\pm\gamma$, with $\gamma\in R$, and $a=1$.
Then $V_{b,a}=\frac{\gamma^2}2$ and thus by unitary equivalence,
$\inf\sigma(H_D)=\inf\sigma(H_S)=\frac{\gamma^2}2$. A direct calculation in this case reveals that $\Omega^+(b,a)=\frac1{4\gamma^2}$; thus,
 $\inf\sigma(H_D)=\frac 1{8\Omega^+(b,a)}$.
On the other hand, consider the case
that $b(x)=-\gamma x$, with $\gamma>0$, and $a=1$.
Then $\lim_{x\to\infty}V_{b,a}(x)=\infty$, so as noted above, $H_S$ and thus also $H_D$ have  compact resolvents.
The unnormalized
 Hermite function $H_1(x)=x$ is an  $L^2$-eigenfunction of $H_D$  corresponding to the eigenvalue
$\gamma$. Since it is positive, it must in fact be the principal eigenvalue. Thus, the bottom
of the spectrum is equal to $\gamma$.
We have
$$
\begin{aligned}
&\Omega^+(b,a)=\sup_{x>0}(\int_0^x\exp(\gamma y^2)dy)(\int_x^\infty\exp(-\gamma y^2)dy)=\\
&\frac1\gamma
\sup_{x>0}(\int_0^x\exp( y^2)dy)(\int_x^\infty\exp(- y^2)dy)\approx\frac{.239}\gamma.
\end{aligned}
$$
Thus, in this case, the bottom of the spectrum is approximately equal to $\frac{.239}{\Omega^+(b,a)}$.
In the case $b(x)=\gamma x$, with $\gamma>0$, and $a=1$, one can check that
the  principal eigenfunction is $x\exp(-\gamma x^2)$, with corresponding principal eigenvalue 2.
One can calculate that $\Omega^+(a,b)\approx\frac{.097}{\gamma}$, and thus the bottom of the spectrum is approximately
equal to $\frac{.194}{\Omega^+(b,a)}$.
Writing the bottom of the spectrum in the form $\frac{C_{b,a}}{\Omega^+(b,a)}$, we
 don't know
 whether the upper bound in the theorem is sharp; namely, $C_{b,a}\le\frac12$.

\medskip

\noindent \bf Remark 2.\rm\  In this remark, we demonstrate how
formulas  \eqref{spec} and \eqref{essspec}  in the case
\eqref{intcondfin} follow from those formulas in the case
\eqref{intcondinfin}, thereby reducing the proof of the theorem to
the case that \eqref{intcondinfin} holds. In the case that
\eqref{intcondfin} holds, define the \it $h$-transform\rm\ of
$H_D$ via the function $h_{b,a}$ in \eqref{hb} by
 $H_D^{h_{b,a}}u=\frac1{h_{b,a}}H_D(h_{b,a}u)$. When written out, one obtains
$H_D^{h_{b,a}}=-\frac12\frac d{dx}a\frac d{dx}-(b+a\frac {h_{b,a}'}{h_{b,a}})\frac d{dx}$.
Letting $B^{h_{b,a}}(x)=\int_0^x(\frac ba+\frac{h_{b,a}'}{h_{b,a}})(y)dy$, one has
$\int^\infty\frac1{a(x)}\exp(-2B^{h_{b,a}}(x))dx=-\int^\infty h_{b,a}^{-2}h_{b,a}'dx=\infty$;
that is, the diffusion coefficient $a$ with the new drift $b+a\frac{h_{b,a}'}{h_{b,a}}$ satisfies \eqref{intcondinfin}.
The spectrum is invariant under $h$-transforms \cite[chapter 4---sections 3 and 10]{P}, so
$\inf\sigma(H_D)=\inf\sigma(H_D^{h_{b,a}})$ and $\inf\sigma_{\text{ess}}(H_D)=\inf\sigma_{\text{ess}}(H_D^{h_{b,a}})$.
These equalities along with the fact that
 \eqref{intcondinfin} holds with the diffusion coefficient $a$ and the drift $(b+a\frac{h_{b,a}'}{h_{b,a}})$ show that one obtains \eqref{spec} and \eqref{essspec} for
 $H_D$ by defining $\Omega^+(b,a)=\Omega^+(b+a\frac{h_{b,a}'}{h_{b,a}},a)$ and
$\hat\Omega^+(b,a)=\hat\Omega^+(b+a\frac{h_{b,a}'}{h_{b,a}},a)$.
 From \eqref{ab}, one has
\begin{equation}
\begin{aligned}
&\Omega^+(b+a\frac{h_{b,a}'}{h_{b,a}},a)=\sup_{x>0}
\left(\int_0^x\frac1{a(y)}\exp(-2B^{h_{b,a}}(y))dy\right)\left(\int_x^\infty\exp(2B^{h_{b,a}}(y))dy\right)\\
&=\sup_{x>0}\left(\int_0^xh_{b,a}^{-2}(y)\frac1{a(y)}\exp(-2B(y))dy\right)\left(\int_x^\infty h^2_{b,a}(y)\exp(2B(y))dy\right);
\end{aligned}
\end{equation}
whence the definition of $\Omega^+(b,a)$ in \eqref{ab-h} in  the case that \eqref{intcondfin} holds, and likewise for $\hat\Omega^+(b,a)$.
\medskip

\noindent \bf Remark 3.\rm\ After finishing this paper, the following related result due to Muckenhoupt \cite{M}, in the context of weighted
Hardy inequalities,
 was brought to our attention.
For $1\le p\le\infty$,
the inequality
\begin{equation}\label{Muck}
\left(\int_0^\infty|U(x)\int_0^xg(t)dt|^pdx\right)^\frac 1p\le C\left( \int_0^\infty |V(x)g(x)|^pdx\right)^\frac1p
\end{equation}
holds for all $g$ and some finite $C$ if and only if
$$
B\equiv\sup_{x>0}\left(\int_x^\infty|U(y)|^pdy\right)^\frac1p\left(\int_0^x|V(y)|^{-p'}dy\right)^\frac1{p'}<\infty,
$$
where $\frac1p+\frac1{p'}=1$, and furthermore, if $C_0$ is the least constant $C$ for which the above inequality holds,
then $B\le C_0\le p^\frac1p (p')^\frac1{p'}B$, for $1<p<\infty$, and $C_0=B$ for $p=1,\infty$.
(The integrals are interpreted according to the usual convention in the case that $p$ or $p'$ is $\infty$.)
Applying this with $p=p'=2$, $U=\exp(B)$ and $V=(\frac a2)^\frac12\exp(B)$,
one concludes that
$\inf\frac{\frac12\int_0^\infty a(f')^2\exp(2B)dx}{\int_0^\infty f^2\exp(2B)dx}$
  lies between
$\left(8\sup_{x>0}(\int_0^x\frac1{a(y)}\exp(-2B(y))dy)(\int_x^\infty\exp(2B(y))dy)\right)^{-1}$
and $\left(2\sup_{x>0}(\int_0^x\frac1{a(y)}\exp(-2B(y))dy)(\int_x^\infty\exp(2B(y))dy)\right)^{-1}$, where
the infimum is over $f\in C^1([0,\infty))$ which satisfy $f(0)=0$.
This is a different variational problem than the one in \eqref{variational} for $\inf\sigma(H_D)$ because the class of admissible
functions here is larger  than in \eqref{variational}. In the case that \eqref{intcondinfin} holds, Theorem \ref{Th} shows that
the same bounds hold for both variational problems, since in this case,
$\Omega^+(b,a)=\sup_{x>0}\left(\int_0^x\frac1{a(y)}\exp(-2B(y))dy\right)\left(\int_x^\infty\exp(2B(y))dy\right)$.
However, when \eqref{intcondfin} holds, $\Omega^+(b,a)$ is defined differently, and the two variational
problems yield different results. Indeed, for example, if $b=1$ and $a=1$, then
one has $\int_x^\infty\exp(2B(y))dy)=\infty$, so the infimum in Muckenhoupt's variational problem is 0; however
by \eqref{ab-h}, one calculates that $\Omega^+(b,a)=1$, and it follows from Theorem \ref{Th} that
the infimum in \eqref{variational} lies between $\frac18$ and $\frac12$. (In fact, in this simple case it can be checked directly
that $\inf\sigma(H_D)=\frac18$.)
Thus,
the integral condition \eqref{intcondfin} turns out to be the lower threshold on the size of
$a\exp(2B)$, the weight that multiplies $(f')^2$ in the variational formulas, so that the two variational formulas,
one over $f\in C^1_0(R^+)$ and one over $f\in C^1([0,\infty))$ satisfying $f(0)=0$, yield different answers.
Muckenhout's proof involves a direct estimation of the integrals in \eqref{Muck}. We prove Theorem \ref{Th}
in a completely different way, as will be seen in sections 3 and 4.

\medskip

\noindent \bf Remark 4.\rm\ Theorem \ref{Th}
 and the reduction noted above in Remark 2 have
 some probabilistic implications, which we now describe.
Let $X(t)$ be generic notation for  a  Markov diffusion
process on the real line. Let
$P_x$ and $E_x$  denote respectively  probabilities and expectations
for the process
 corresponding to the
operator $-H_D=\frac12\frac d{dx}a\frac d{dx}+b\frac d{dx}$ on $(0,\infty)$, starting at $x>0$  and killed at time
\begin{equation}\label{hitting}
\tau_0\equiv\inf\{t\ge0:X(t)=0\},
\end{equation}
the first hitting time of 0. Then
$P_x(\tau_0<\infty)=1$, for $x>0$,  if and only if \eqref{intcondinfin} holds \cite[chapter 5]{P}.

Consider first the case that \eqref{intcondinfin} holds.
At the end of  section 3 we show that
\begin{equation}\label{exphitting}
\inf\sigma(H_D)=\sup\{\lambda\ge0: E_x\exp(\lambda\tau_0)<\infty\}, \ x>0.
\end{equation}
 Thus, in the case that \eqref{intcondinfin} holds,  \eqref{spec}  gives an  explicit formula up to the multiplicative
 constant 4 for
$\sup\{\lambda\ge0: E_x\exp(\lambda\tau_0)<\infty\}$.

Now consider the case that \eqref{intcondfin} holds.
In this case, $P_x(\tau_0<\infty)=\frac{h_{b,a}(x)}{h_{b,a}(0)}$ \cite[chapter 5]{P}. The original process, conditioned on
$\{\tau_0<\infty\}$, is itself a Markov diffusion process and it corresponds to the
$h$-transformed operator  $-H_D^{h_{b,a}}$ defined in Remark 2 \cite[chapter 7]{P}. Let $E_x^{h_{b,a}}$ denote expectations for this conditioned process
starting from $x>0$.
Then it follows from Remark 2 and \eqref{exphitting} that
in the case that  \eqref{intcondfin} holds, one has
\begin{equation}\label{exphittingh}
\inf\sigma(H_D)=\sup\{\lambda\ge0: E^{h_{b,a}}_x\exp(\lambda\tau_0)<\infty\},\ x>0.
\end{equation}

Note from \eqref{ab} and  \eqref{spec} that when \eqref{intcondinfin} holds, a necessary condition for $\inf\sigma(H_D)>0$ is that
$\int^\infty\exp(2B(y))dy<\infty$. This integral condition is equivalent to $E_x\tau_0<\infty$, for $x>0$
\cite[chapter 5---section 1]{P}.
Thus, when $P_x(\tau_0<\infty)=1$  holds, the finiteness of $E_x\tau_0$  is a necessary
condition (but not a sufficient one) for   $\inf\sigma(H_D)>0$.
Similarly,  when \eqref{intcondfin} holds (in which case
$P^{h_{b,a}}(\tau_0<\infty)=1$), the finiteness of $E_x^{h_{b,a}}\tau_0$ is a necessary condition (but not a sufficient one) for
$\inf\sigma(H_D)>0$. (Of course, this can also be seen from \eqref{exphitting} and \eqref{exphittingh}---if the first moment
does not exist, then a fortiori no exponential moment exists.)

An alternative probabilistic representation of $\inf\sigma(H_D)$ is this:
\begin{equation}\label{probrep}
\inf\sigma(H_D)=-\lim_{n\to\infty}\lim_{t\to\infty}\frac1t\log P_x(\tau_0\wedge\tau_n>t),\ x>0,
\end{equation}
where $\tau_n=\inf\{t\ge0:X(t)=n\}$. (This formula can  be found essentially in \cite[chapter 4]{P}.)

Formulas \eqref{exphitting} and \eqref{exphittingh} give a probabilistic representation for the bottom
of the spectrum of $H_D$. One can also give a similar probabilistic representation for the bottom of the essential
spectrum.
It follows from \eqref{exphitting} and \eqref{limitalpha} in section 3 that if \eqref{intcondinfin} holds, then
\begin{equation*}
\inf\sigma_{\text{ess}}(H_D)=\lim_{l\to\infty}(\sup\{\lambda\ge0: E_x\exp(\lambda\tau_l)<\infty\ \ x>l\}),
\end{equation*}
while if \eqref{intcondfin} holds, then
\begin{equation*}
\inf\sigma_{\text{ess}}(H_D)=\lim_{l\to\infty}(\sup\{\lambda\ge0: E^{h_{b,a}}_x\exp(\lambda\tau_l)<\infty,\  \ x>l\}).
\end{equation*}

\medskip

\noindent \bf Remark 5.\rm\ It follows from the theorem that
$\inf\sigma(H_D)$ and $\inf\sigma_{\text{ess}}(H_D)$ depend on $a$ and $b$
only through $a$ and $B$.

\medskip

\noindent \bf Remark 6.\rm\ For the duration of this remark, we consider $a$ to be fixed.
 By a standard comparison theorem for diffusions,
it follows that
 over the class of drifts $b$ satisfying
\eqref{intcondinfin}, the distribution of $\tau_0$ is stochastically increasing with $b$. Thus, from \eqref{exphitting}, it follows that
 $\inf\sigma(H_D)$ and $\inf\sigma_{\text{ess}}(H_D)$ are nonincreasing over the class of drifts $b$ satisfying \eqref{intcondinfin}.
That is,  over drifts satisfying \eqref{intcondinfin}, the more inward toward 0 the drift, the larger the bottom of the spectrum
and the bottom of the essential spectrum.
(It is not hard to verify that the function
$\left(\int_0^x\frac1{a(y)}\exp(-2B(y))dy\right)\left(\int_x^\infty\exp(2B(y))dy\right)$ appearing in the definition of $\Omega^+(b,a)$
 is nondecreasing in
$b$ over the class of drifts $b$ satisfying \eqref{intcondinfin}, but this is not quite enough to arrive at the result in the above sentence.)
Despite the above fact and despite Remark 5,
 it is \it not\rm\ true that  $\inf\sigma(H_D)$ and $\inf\sigma_{\text{ess}}(H_D)$ are nonincreasing as functions of $B$ over the class
of $b$ satisfying \eqref{intcondinfin}. An example will  be given at the end of  section 2.

We don't know whether   $\inf\sigma(H_D)$ and $\inf\sigma_{\text{ess}}(H_D)$ are nondecreasing over the entire class of drifts $b$ satisfying \eqref{intcondfin},
so that the more outward toward infinity the drift, the larger
the bottom of the spectrum
and the bottom of the essential spectrum.
To prove that this is true, it would suffice
to show that $b+a\frac{h'_{b,a}}{h_{b,a}}$ is nonincreasing in $b$ over the class of drifts
satisfying \eqref{intcondfin}---that this would suffice follows from
\eqref{exphittingh} and the argument above for the class of drifts satisfying \eqref{intcondinfin}.
What is known is this \cite{P93}:
\begin{equation}\label{opposite}
\begin{aligned}
&\text{For a wide class of}\ a\ \text{ and}\ b\ \text{ which satisfy \eqref{intcondfin} and for which}\ b\ \text{ is on a}\\
 &\text{larger
order than}\  \frac{a(x)}x,\ \text{ one has}\
b+a\frac{h'_{b,a}}{h_{b,a}}=-b+O(\frac{a(x)}x)\ \text{ as}\ x\to\infty.
\end{aligned}
\end{equation}
This formula will be useful for one of the calculations in section 2.
\bigskip

We now turn to the case of the whole line. Let $0<a\in C^1(R)$ and  $b\in C(R)$, and define $B(x)=\int_0^x\frac ba(y)dy$.
Let
 $H_D=-\frac12\frac d{dx}a\frac d{dx}-b\frac d{dx}$  and  consider the self-adjoint
realization on $L^2(R,\exp(2B)dx)$ obtained via the Friedrichs extension of the closure of
the quadratic form $Q_D(f,g)=\frac12\int_{-\infty}^\infty (f'ag')\exp(2B)dx$, for $f,g\in C_0^1(R)$.
In the case that $b\in C^1(R)$, $\frac{b^2}a+b'$ is bounded from below and $\int^\infty a(x)dx=\int_{-\infty}a(x)dx=\infty$, define
$H_S=-\frac12\frac d{dx}a\frac d{dx}+V_{b,a}$ to be
the self-adjoint operator obtained via the  Friedrichs extension of the
closure of the
quadratic form $Q_S(f,g)=\frac12\int_{-\infty}^\infty (f'ag')dx+
\int_{-\infty}^\infty V_{b,a}fgdx$,
where $V_{b,a}=\frac12(\frac{b^2}a+b')$ and $f,g\in C_0^1(R)$.

The first of the two theorems below treats $\inf\sigma_{\text{ess}}(H_D)$ and the second one treats $\inf\sigma(H_D)$.
The proofs of these results will be derived in just a few lines from the proof of Theorem \ref{Th}.

\begin{theorem}\label{Th2}
Let $a\in C^1(R)$ and
$b\in C(R)$. Define
$$
B(x)=\int_0^x\frac ba(y)dy.
$$
Consider the  self-adjoint diffusion operator
$$
H_D=-\frac12\frac d{dx}a\frac d{dx}-b\frac d{dx}=
-\frac12\exp(-2B)\frac d{dx}a\exp(2B)\frac d{dx}
$$
 on
$L^2(R,\exp(2B)dx)$.

\noindent If
$b\in C^1(R)$,
$\frac{b^2}a+b'$ is bounded from below and $\int^\infty a(x)dx=\int_{-\infty} a(x)dx=\infty$, consider also  the self-adjoint
Schr\"odinger-type operator
$$
H_S=-\frac12\frac d{dx}a\frac d{dx}+\frac12(\frac{b^2}a+b')
$$
on $L^2(R)$.
Let $\Omega^+(b,a)$ be as in Theorem \ref{Th} and define $\Omega^-(b,a)$ in exactly the same way, using the half-line $(-\infty,0)$ instead
of $(0,\infty)$.
Let
\begin{equation*}
\hat \Omega(b,a)=\max(\hat \Omega^+(b,a),\hat \Omega^-(b,a)).
\end{equation*}
Then
\begin{equation*}
\frac1{8\hat \Omega(b,a)}\le\inf\sigma_{\text{ess}}(H_D)=\inf\sigma_{\text{ess}}(H_S)\le\frac1{2\hat \Omega(b,a)}.
\end{equation*}
In particular, $H_D$ and $H_S$ possess compact
resolvents if and only if $\hat \Omega^+(b,a)=\hat \Omega^-(b,a)=0$.
\end{theorem}
\noindent \bf Remark 7.\rm\
The diffusion is
positive recurrent  if and only if $\int_R\exp(2B(x))dx<\infty$ \cite[chapter 5]{P}. It follows from Theorem 2 that $\inf\sigma_{\text{ess}}(H_B)=0$
if the diffusion is not positive recurrent. (See also the third to the last paragraph of Remark 4.)

\medskip

\begin{theorem}\label{Th3}
Let $a\in C^1(R)$ and
$b\in C(R)$. Define
$$
B(x)=\int_0^x\frac ba(y)dy.
$$
Consider the  self-adjoint diffusion operator
$$
H_D=-\frac12\frac d{dx}a\frac d{dx}-b\frac d{dx}=
-\frac12\exp(-2B)\frac d{dx}a\exp(2B)\frac d{dx}
$$
 on
$L^2(R,\exp(2B)dx)$.

\noindent If
$b\in C^1(R)$,
$\frac{b^2}a+b'$ is bounded from below and $\int^\infty a(x)dx=\int_{-\infty} a(x)dx=\infty$, consider also  the self-adjoint
Schr\"odinger-type operator
$$
H_S=-\frac12\frac d{dx}a\frac d{dx}+\frac12(\frac{b^2}a+b')
$$
on $L^2(R)$.

\noindent If
\begin{equation}\label{recurrent}
\int^\infty\frac1{a(x)}\exp(-2B(x))dx=\int_{-\infty}\frac1{a(x)}\exp(-2B(x))dx=\infty,
\end{equation}
define
\begin{equation}
\Omega(b,a)=\infty.
\end{equation}
If
\begin{equation}\label{trans-}
\int^\infty\frac1{a(x)}\exp(-2B(x))dx=\infty\ \ \text{and}\ \ \int_{-\infty}\frac1{a(x)}\exp(-2B(x))dx<\infty,
\end{equation}
define
\begin{equation*}
\Omega(b,a)=\sup_{x\in R}\left(\int_{-\infty}^x\frac1{a(y)}\exp(-2B(y))dy\right)\left(\int_x^\infty\exp(2B(y))dy\right).
\end{equation*}
If
\begin{equation}\label{trans+}
\int^\infty\frac1{a(x)}\exp(-2B(x))dx<\infty\ \ \text{and}\ \ \int_{-\infty}\frac1{a(x)}\exp(-2B(x))dx=\infty,
\end{equation}
define
\begin{equation*}
\Omega(b,a)=\sup_{x\in R}\left(\int_x^\infty\frac1{a(y)}\exp(-2B(y))dy\right)\left(\int_{-\infty}^x\exp(2B(y))dy\right).
\end{equation*}
If
\begin{equation}\label{trans+-}
\int^\infty\frac1{a(x)}\exp(-2B(x))dx<\infty\ \  \text{and}\ \ \int_{-\infty}\frac1{a(x)}\exp(-2B(x))dx<\infty,
\end{equation}
let
\begin{equation*}
h_{b,a}(x)=\int_x^\infty\frac1{a(y)}\exp(-2B(y))dy
\end{equation*}
and define
\begin{equation*}
\begin{aligned}
&\Omega(b,a)=\sup_{x\in R}\left(\int_{-\infty}^xh^{-2}_{b,a}(y)\frac1{a(y)}\exp(-2B(y))dy\right)\left(\int_x^\infty h^2_{b,a}(y)\exp(2B(y))dy\right)\\
&=\sup_{x\in R}\left(h_{b,a}^{-1}(x)-h_{b,a}^{-1}(-\infty)\right)\left(\int_x^\infty h^2_{b,a}(y)\exp(2B(y))dy\right).
\end{aligned}
\end{equation*}
 Then
\begin{equation*}
\frac1{8\Omega(b,a)}\le\inf\sigma(H_D)=\inf\sigma(H_S)\le\frac1{2\Omega(b,a)}.
\end{equation*}

\end{theorem}

\noindent \bf Remark 8.\rm\ The diffusion process $X(t)$ corresponding to $-H_D$ is recurrent if \eqref{recurrent}
holds and is transient otherwise. In the transient case, if \eqref{trans-} holds, then $P_x(\lim_{t\to\infty}X(t)=-\infty)=1$;
if \eqref{trans+} holds, then $P_x(\lim_{t\to\infty}X(t)=\infty)=1$;
if \eqref{trans+-} holds, then
$P_x(\lim_{t\to\infty}X(t)=-\infty)=1-P_x(\lim_{t\to\infty}X(t)=\infty)=\frac{h_{b,a}(x)}{h_{b,a}(-\infty)}$.
(For these results,
see \cite[chapter 5]{P}.)
It follows from Theorem \ref{Th3} that $\inf\sigma(H_D)=0$ if the diffusion is recurrent.

\medskip

\noindent \bf Remark 9.\rm\
Similar to \eqref{probrep}, one has the following probabilistic representation of $\inf\sigma(H_D)$:
\begin{equation}\label{probrep2}
\inf\sigma(H_D)=-\lim_{n\to\infty}\lim_{t\to\infty}\frac1t\log P_x(\tau_{-n}\wedge\tau_n>t),\ x\in R.
\end{equation}
\medskip

In section 2 we give some applications of Theorems \ref{Th}-\ref{Th3}.
In section 3 we prove Theorem \ref{Th}, postponing the proof of a key proposition to section 4.
After the proof of Theorem \ref{Th}  we give the quick proofs of Theorems \ref{Th2} and \ref{Th3}.
We also prove \eqref{exphitting} in section 3.
Finally, in section 5 we show how the one-dimensional result can be used to obtain   spectral estimates
 for self-adjoint, multi-dimensional  diffusion operators

\section{Examples}

\noindent \it The Bottom of the Spectrum.\rm\
One can use Theorem \ref{Th} to study the way $\inf\sigma(H_D)$ scales in the parameters $\gamma$ and $\nu$ when
$b$ is of the form $b=\gamma b_0$ and $a$ is of the form $a=\nu a_0$.
We first consider the effect of the drift alone.
Consider for example the following
two cases on $R^+$ or on $R$:
\begin{equation}\label{balone}
b(x)=-\gamma(1+|x|)^l\ \text{and}\ a(x)=1, \ \gamma>0,\ l\in R,
\end{equation}
\begin{equation}\label{baloneagain}
b(x)=-\gamma |x|^l\ \text{and}\ a(x)=1, \ \gamma>0,\ l\ge0.
\end{equation}

\begin{proposition}\label{justb}
Consider $H_D$ on $R^+$ or on $R$.

\noindent 1. Assume that \eqref{balone} holds.

i. If $l<0$, then $\inf\sigma(H_D)=0$;

ii. If $l\ge0$, then there exist constants $c_l, C_l>0$ such that
\begin{equation*}
c_l\gamma^2\le\inf\sigma(H_D)\le C_l\gamma^2,\ \gamma>1
\end{equation*}
and
\begin{equation*}
c_l\gamma^\frac2{l+1}\le\inf\sigma(H_D)\le C_l\gamma^\frac2{l+1},\ 0<\gamma\le1.
\end{equation*}

\noindent 2. Assume that \eqref{baloneagain} holds.
Then
\begin{equation*}
\frac1{8C_l}\gamma^\frac2{1+l}\le\inf\sigma(H_D)\le \frac1{2C_l}\gamma^\frac2{1+l},\ \gamma>0,
\end{equation*}
where
$$
C_l=
\begin{cases}&\sup_{x>0}\left(\int_0^x\exp(\frac{2z^{l+1}}{l+1})dz\right)\left(\int_x^\infty\exp(-\frac{2z^{l+1}}{l+1})dz\right)
\ \text{on}\ R^+;\\
&C_l=\left(\int_0^\infty\exp(-\frac{2z^{l+1}}{l+1})dz\right)^2 \ \text{on}\ R.
\end{cases}
$$
\end{proposition}

\noindent \bf Remark 10.\rm\
Note that both on $R^+$ and on $R$, the rate of growth of $\inf\sigma(H_D)$ for large $\gamma$ is on a slower order
for the drift in \eqref{baloneagain}  than for the drift in \eqref{balone}.
The probabilistic explanation for this follows from the formulas \eqref{probrep} and \eqref{probrep2}
and the fact that the latter drifts are small in a ($\gamma$-dependent) neighborhood of 0, even as $\gamma$ becomes large.
Note also that for the drift in \eqref{balone}, the scaling power is different for $\gamma\ll1$ than for $\gamma\gg1$.

 The bounds on the infimum of the spectrum in Proposition \ref{justb} also hold for the corresponding Schr\"odinger operator,
$H_S=-\frac12\frac{d^2}{dx^2}+V$, where $V=\frac12\gamma^2(1+x)^{2l}-\frac12\gamma l(1+x)^{l-1}$ in the case of \eqref{balone} on $R^+$
and $V(x)=\frac12\gamma^2x^{2l}-\frac12\gamma x^{l-1}$
 in the case of \eqref{baloneagain} on $R^+$, and  similarly for $R$.

\medskip

We now consider simultaneous scaling in $a$ and $b$.
Consider the following case on $R^+$ and on $R$:
\begin{equation}\label{ab-bottom}
\begin{aligned}
&b(x)=-\gamma(1+|x|)^l\ \ \text{and}\ a(x)=\nu(1+|x|)^k,\\
&\text{where} \ \gamma,\nu>0,\ l,k\in R,\
\text{with}\ l-k>-1\ \text{and}\ 2l-k\ge0.
\end{aligned}
\end{equation}
(Note that when $\nu=1$ and $k=0$, \eqref{ab-bottom} reduces to \eqref{balone} with $l\ge0$.)
If $2l-k<0$ or if $l-k<-1$, then one can show that $\inf\sigma(H_D)=0$.
\begin{proposition}\label{scaling}
Consider $H_D$ on $R^+$ or on $R$.
 Assume that  \eqref{ab-bottom} holds.
 There exist constants  $c_{l,k}, C_{l,k}>0$ such that
$$
c_{l,k}\frac{\gamma^2}\nu\le\inf\sigma(H_D)\le C_{l,k}\frac{\gamma^2}\nu, \ 0<\nu<\gamma,
$$
and
$$
c_{l,k}(\frac{\gamma^{2-k}}{\nu^{1-l}})^\frac1{l-k+1}\le
\inf\sigma(H_D)\le C_{l,k}(\frac{\gamma^{2-k}}{\nu^{1-l}})^\frac1{l-k+1},\ 0<\gamma\le\nu.
$$

\end{proposition}

\noindent \bf Remark 11.\rm\ Note that when $\gamma\le\nu$,  the scaling dependence on the coefficient $\gamma$ of the drift $b$ has three dramatically different phases, depending on
whether the exponent $k$ of the diffusion coefficient satisfies $k<2$, $k=2$ or $k>2$, while the scaling dependence of the coefficient
$\nu$ of the diffusion coefficient has three dramatically different phases, depending on whether the exponent $l$ of the drift satisfies
$l<1$, $l=1$ or $l>1$.  However, when $\gamma>\nu$, there is
only one scaling phase, and it is independent of the exponents $l$ and $k$.

The bounds on the infimum of the spectrum in Proposition \ref{scaling} also hold for the corresponding Schr\"odinger-type operator,
$H_S=-\frac12\frac d{dx}(\nu(1+|x|)^k)\frac d{dx}+V$, where $V=\frac12\frac{\gamma^2}\nu(1+x)^{2l-k}-\frac12\gamma l(1+x)^{l-1}$ in the case of $R^+$, and similarly for $R$.
The parameter dependence in Proposition \ref{scaling}
does not seem at all apparent from looking at this operator.

\medskip
 We give the proof of Proposition \ref{justb}; the proof of Proposition \ref{scaling} is similar.

\noindent \bf Proof of Proposition \ref{justb}.\rm\
We prove the proposition in the case of $R^+$; the case of $R$ is handled similarly.
To prove part 2, one simply makes an appropriate change of variables in the formula for $\Omega^+(-\gamma x^l,1)$ and applies Theorem \ref{Th}.
To get the explicit form of $C_l$ in the case of $R$, one needs to do a little bit more analysis to show that the
supremum over $x\in R$ occurs at $x=0$.

We now prove part 1.
If $l\neq-1$, then
\begin{equation}\label{omegaex}
\begin{aligned}
&\Omega^+(-\gamma(1+x)^l,1)=\\
&\sup_{x>0}\left(\int_0^x\exp(2\gamma\frac{(1+y)^{l+1}}{l+1})dy\right)
\left(\int_x^\infty\exp(-2\gamma\frac{(1+y)^{l+1}}{l+1})dy\right).
\end{aligned}
\end{equation}
For $l<-1$ the right hand integral is $\infty$ so $\Omega^+(-\gamma(1+x)^l,1)=\infty$. Now consider $-1<l<0$.
Applying L'H\^opital's rule to the quotients
$$
\frac{\int_0^x\exp(2\gamma\frac{(1+y)^{l+1}}{l+1})dy}{(1+x)^{-l}\exp(2\gamma\frac{(1+x)^{l+1}}{l+1})},\
\frac{\int_x^\infty\exp(-2\gamma\frac{(1+y)^{l+1}}{l+1})dy}{(1+x)^{-l}\exp(-2\gamma\frac{(1+x)^{l+1}}{l+1})},\
$$
shows that $\int_0^x\exp(2\gamma\frac{(1+y)^{l+1}}{l+1})dy\sim (2\gamma)^{-1}(1+x)^{-l}\exp(2\gamma\frac{(1+x)^{l+1}}{l+1})$ and \newline
$\int_x^\infty\exp(-2\gamma\frac{(1+y)^{l+1}}{l+1})dy\sim (2\gamma)^{-1}(1+x)^{-l}\exp(-2\gamma\frac{(1+x)^{l+1}}{l+1})$, as
$x\to\infty$.
This shows that the supremum  in \eqref{omegaex} is $\infty$; thus $\Omega^+(-\gamma(1+x)^l,1)=\infty$. One obtains $\Omega^+(-\gamma(1+x)^l,1)=\infty$ similarly
in the case $l=-1$. Applying Theorem \ref{Th} now completes the proof of part 1-i.

Consider now part 1-ii; that is, the case $l\ge0$. Making the change of variables $z=\gamma^{\frac1{l+1}}(1+y)$, one obtains from \eqref{omegaex},
\begin{equation}\label{omegaex2}
\begin{aligned}
\Omega^+(-\gamma(1+x)^l,1)=\gamma^{-\frac2{l+1}}\sup_{x>\gamma^{\frac1{l+1}}}\left(\int^x_{\gamma^{\frac1{l+1}}}\exp(\frac{z^{l+1}}{l+1})dz\right)
\left(\int_x^\infty \exp(-\frac{z^{l+1}}{l+1})dz\right).
\end{aligned}
\end{equation}
If $l=0$, the integrals on the right hand side of \eqref{omegaex2} can be calculated explicitly. One finds that the supremum above is equal to 1. Part
1-ii in the case $l=0$ now follows from \eqref{omegaex2} and Thereom \ref{Th}.
From now on, we assume that $l>0$.
Applying L'H\^opital's rule in the manner noted above shows that
\begin{equation}\label{lhopital}
\begin{aligned}
&\int_x^\infty\exp(-\frac{z^{l+1}}{l+1})dz\sim x^{-l}\exp(-\frac{x^{l+1}}{l+1}),\ \ \text{as}\ x\to\infty;\\
&\int_0^x\exp(-\frac{z^{l+1}}{l+1})dz\sim x^{-l}\exp(\frac{x^{l+1}}{l+1}),\ \ \text{as}\ x\to\infty.
\end{aligned}
\end{equation}
From \eqref{lhopital} it follows that there exist constants $d_l,D_l>0$ such that
\begin{equation}\label{omegabounds}
d_l\le\sup_{x>\gamma^{\frac1{l+1}}}\left(\int^x_{\gamma^{\frac1{l+1}}}\exp(\frac{z^{l+1}}{l+1})dz\right)
\left(\int_x^\infty \exp(-\frac{z^{l+1}}{l+1})dz\right)\le D_l, \ 0<\gamma\le1.
\end{equation}
Part 1-ii in the case that $0<\gamma\le 1$ now follows from \eqref{omegaex2}, \eqref{omegabounds} and Theorem \ref{Th}.

Now consider part 1-ii in the case that $\gamma>1$. Clearly,
\begin{equation}\label{twosided}
\begin{aligned}
&\left(\int^{2\gamma^{\frac1{l+1}}}_{\gamma^{\frac1{l+1}}}\exp(\frac{z^{l+1}}{l+1})dz\right)
\left(\int_{2\gamma^{\frac1{l+1}}}^\infty \exp(-\frac{z^{l+1}}{l+1})dz\right)\le\\
&\sup_{x>\gamma^{\frac1{l+1}}}\left(\int^x_{\gamma^{\frac1{l+1}}}\exp(\frac{z^{l+1}}{l+1})dz\right)
\left(\int_x^\infty \exp(-\frac{z^{l+1}}{l+1})dz\right)\le\\
&\sup_{x>\gamma^{\frac1{l+1}}}\left(\int^x_0\exp(\frac{z^{l+1}}{l+1})dz\right)
\left(\int_x^\infty \exp(-\frac{z^{l+1}}{l+1})dz\right).
\end{aligned}
\end{equation}
Using \eqref{lhopital} to estimate the left most and right most terms in \eqref{twosided},  it follows that there exist constants $d_l,D_l>0$ such that
\begin{equation}\label{omegaboundsagain}
\begin{aligned}
&d_l\gamma^{-\frac{2l}{l+1}}\le\sup_{x>\gamma^{\frac1{l+1}}}\left(\int^x_{\gamma^{\frac1{l+1}}}\exp(\frac{z^{l+1}}{l+1})dz\right)
\left(\int_x^\infty \exp(-\frac{z^{l+1}}{l+1})dz\right)\le D_l\gamma^{-\frac{2l}{l+1}},\\
&\text{for}\ \gamma>1.
\end{aligned}
\end{equation}
Part 1-ii in the case $\gamma>1$ now follows from \eqref{omegaex2}, \eqref{omegaboundsagain} and Theorem \ref{Th}.
\hfill$\square$

\medskip

Theorem \ref{Th}  allows one to compute the bottom of the spectrum \it exactly\rm\ for an ad hoc class of Schr\"odinger operators,
$H=-\frac12\frac{d^2}{dx^2}+V$.
Indeed, it follows from the theorem that if $a$ and $b$ satisfy \eqref{intcondinfin} and $\int^\infty\exp(2B)dx=\infty$, then
$\inf\sigma(H_D)=\inf\sigma(H_S)=0$.
Let $u=\exp(g)$, where $g$ is bounded, and define $V=\frac{u''}{2u}=\frac12((g')^2+g'')$,  $b=\frac{u'}u=g'$ and $a=1$.
 Then $b$ satisfies the above conditions and $H_S=-\frac12\frac{d^2}{dx^2}+V$.
  Thus,
$\inf\sigma(-\frac12\frac{d^2}{dx^2}+\frac12((g')^2+g''))=0$, for all bounded $g$.
\it In particular,  if $g$ is periodic and not constant, then $\lim_{x\to\infty}\frac1x\int_0^xV(y)dy>0$ but the bottom of the
spectrum is 0.\rm\

Note that either $\hat \Omega^+(b,a)=\Omega^+(b,a)=\infty$, or $\hat \Omega^+(b,a),\Omega^+(b,a)<\infty$; thus, in $R^+$ either both the bottom of the spectrum and bottom
of the essential
spectrum equal 0,
or else neither of them does.
It is not hard to construct examples where the bottom of the spectrum and the bottom of the essential spectrum  are both positive and finite
but don't coincide.
For example, let $a=1$ and let $b(x)=-1$, for $x\ge3$. Since $\hat \Omega^+(b,1)$ does not depend on $\{b(x):0\le x\le 3\}$, we
have   $\hat \Omega^+(b,1)=\frac14$. Let $b(x)=-n$, for $1\le x\le2$, and $b(x)\le-1$ everywhere. Then the term $\int_0^3\exp(-2B(y))dy$ can
be made arbitrarily large by choosing $n$ arbitrarily large, and thus for sufficiently large $n$,
$$
\Omega^+(b,1)=\sup_{x>0}\left(\int_0^x\exp(-2B(y))dy\right)\left(\int_x^\infty\exp(2B(y))dy\right)>\frac14=\hat \Omega^+(b,1).
$$
\medskip

\noindent \it The Bottom of the Essential Spectrum.\rm\
We consider operators on  $R^+$. The examples can easily be extended to operators on $R$ by making the analysis on $R^+$
and on $R^-$ separately, and applying Theorem \ref{Th2}.
Consider first the case that
\begin{equation}\label{example-ab}
\begin{aligned}
&b(x)=-\gamma(1+x)^l\ \ \text{and}\ a(x)=\nu(1+x)^k,\ \ \gamma,\nu>0,\ l,k\in R,\\
&\text{with}\ l-k>-1,\ \text{or}\ l-k=-1 \ \text{and}\ k\le1+ \frac{2\gamma}\nu,\ \text{or}\
l-k<-1\ \text{and}\ k\le 1.
\end{aligned}
\end{equation}
The set of possible conditions on $l,k$ above  are exactly those for which   \eqref{intcondinfin}  holds.
One can obtain the asymptotic behavior of $\int_x^\infty\exp(2B(y))dy$ and of
$\int_0^x\frac1{a(y)}\exp(-2B(y))dy$ by
applying L'H\^opital's rule respectively to \newline $\frac{\int_x^\infty\exp(2B(y))dy}{b^{-1}(x)\exp(2B(x))}$
and $\frac{\int_0^x\frac1{a(y)}\exp(-2B(y))dy}{b^{-1}(x)\frac1{a(x)}\exp(-2B(x))}$.
Calculating and applying Theorem \ref{Th}, one obtains the following result.
\medskip

\begin{proposition}\label{exab}
Consider $H_D$ on $R^+$.
Let $a$ and $b$ satisfy \eqref{example-ab}.

\noindent 1. Assume that $l-k<-1$ or that $l-k=-1$ and $\frac\gamma\nu\le\frac12$. Then $\inf\sigma_{\text{ess}}(H_D)=0$.

\noindent 2. Assume that $l-k=-1$ and $\frac\gamma\nu>\frac12$.

i. If $k>2$, then $\sigma_{\text{ess}}(H_D)=\emptyset$;

ii. If $k=2$, then $0<\inf\sigma_{\text{ess}}(H_D)<\infty$;

iii. If $k<2$, then $\inf\sigma_{\text{ess}}(H_D)=0$.

\noindent 3. Assume that $l-k>-1$.

i. If $2l-k>0$, then $\sigma_{\text{ess}}(H_D)=\emptyset$;

ii. If $2l-k=0$, then $0<\inf\sigma_{\text{ess}}(H_D)<\infty$;

iii. If $2l-k<0$, then $\inf\sigma_{\text{ess}}(H_D)=0$.

In particular, $H_D$ possesses a compact resolvent if and only if 2-i or 3-i holds.
\end{proposition}

The bounds on the infimum of the essential spectrum in Proposition \ref{exab} also hold for the corresponding Schr\"odinger-type
operator $H_S=-\frac12\frac d{dx}(\nu(1+x)^k)\frac d{dx}+V$, where $V=\frac12\frac{\gamma^2}\nu(1+x)^{2l-k}-\frac12\gamma lx^{l-1}$.
For certain values of the parameters, the results in Proposition \ref{exab} can be deduced directly from looking at $H_S$.
For example, if $k=0$ and $l<1$, then $\lim_{x\to\infty}V(x)$ equals $\infty$ if $l>0$ and is equal to $\frac12\frac{\gamma^2}\nu$ if $l=0$.
It follows from standard perturbations results, mentioned in the first section, that in the former case
$\sigma_{\text{ess}}(H_S)=\emptyset$ and in the latter case $\inf\sigma_{\text{ess}}(H_S)=\frac12\frac{\gamma^2}\nu$.

However, in fact, Theorem \ref{Th}
 allows one to come to the same type of conclusions as in Proposition \ref{exab} in the case that $a$ and $b$ satisfy one of the following
 general conditions:
\begin{equation}\label{example-abgen}
\begin{aligned}
 &c_1(1+x)^k\le a(x)\le c_2(1+x)^k,\  \ k\in R\\
 &   -c_2(1+x)^m\le \int_0^x\frac{b(y)}{(1+y)^k}dy\le -c_1(1+x)^m, \ \text{for large}\ x,\ \  m>0, \ 0<c_1<c_2;\\
 \end{aligned}
 \end{equation}
or
\begin{equation}\label{example-abgen0}
\begin{aligned}
 &c_1(1+x)^k\le a(x)\le c_2(1+x)^k,\ k\le1\\
    &  \int_0^x\frac{b(y)}{(1+y)^k}dy\ \text{is bounded in}\ x.
 \end{aligned}
 \end{equation}
It is easy to check that under \eqref{example-abgen} or \eqref{example-abgen0}, $a$ and $b$ satisfy \eqref{intcondinfin}.

  Note that now $b$ can be locally erratic, and the bottom of the essential spectrum cannot be deduced directly by looking at $H_S$.

\begin{proposition}\label{exabgen}
Consider $H_D$ on $R^+$.

 \noindent 1. Assume that   $a$ and $b$ satisfy \eqref{example-abgen}.

i. If $2m+k-2>0$,    then $\sigma_{\text{ess}}(H_D)=\emptyset$;

ii. If $2m+k-2=0$, then $0<\inf\sigma_{\text{ess}}(H_D)<\infty$;

iii. If $2m+k-2<0$, then
$\inf\sigma_{\text{ess}}(H_D)=0$.

\noindent 2. Assume that $a$ and $b$ satisfy \eqref{example-abgen0}. Then $\inf\sigma_{\text{ess}}(H_D)=0$.

In particular, $H_D$ possesses a compact resolvent if and only if 1-i holds.

\end{proposition}
To prove Proposition \ref{exabgen}, one makes  the same kind of analysis used for
the proof of Proposition \ref{exab},
  along with the  following monotonicity property which is easy to verify:
for fixed $a$, if \eqref{intcondinfin} holds, then for any $x_0>0$,
$\hat \Omega^+(b,a)$ does not depend on $\{b(x), 0\le x\le x_0\}$ and  it is nondecreasing as a function of $\{b(x), x>x_0\}$.

In Propositions \ref{exab} and \ref{exabgen}, the coefficients $a$ and $b$ are such that
\eqref{intcondinfin} holds.
When \eqref{intcondfin} holds instead, the analysis is more complicated. We state the following analogous
result for the case that \eqref{intcondfin} holds. Consider the following
analog of \eqref{example-abgen}:
\begin{equation}\label{example-abgenout}
\begin{aligned}
 &c_1(1+x)^k\le a(x)\le c_2(1+x)^k,\  \ k\in R\\
 &   c_1(1+x)^m\le \int_0^x\frac{b(y)}{(1+y)^k}dy\le c_2(1+x)^m, \ \text{for large}\ x,\ \  m>0, \ 0<c_1<c_2,\\
 \end{aligned}
 \end{equation}
 and the following analog of \eqref{example-abgen0}:
\begin{equation}\label{example-abgen0out}
\begin{aligned}
 &c_1(1+x)^k\le a(x)\le c_2(1+x)^k,\ k>1\\
    &  \int_0^x\frac{b(y)}{(1+y)^k}dy\ \text{is bounded in}\ x.
 \end{aligned}
 \end{equation}
It can be checked that under \eqref{example-abgenout} or \eqref{example-abgen0out}, $a$ and $b$ satisfy \eqref{intcondfin}.

\begin{proposition}\label{exab-out}
Consider $H_D$ on $R^+$.
Under some mild regularity conditions on $a$ and $b$ one has the following:

\noindent 1. Assume that   $a$ and $b$ satisfy \eqref{example-abgenout}.

i. If $2m+k-2>0$,    then $\sigma_{\text{ess}}(H_D)=\emptyset$;

ii. If $2m+k-2=0$, then $0<\inf\sigma_{\text{ess}}(H_D)<\infty$;

iii. If $2m+k-2<0$, then
$\inf\sigma_{\text{ess}}(H_D)=0$.

\noindent 2. Assume that $a$ and $b$ satisfy \eqref{example-abgen0out}. Then $\inf\sigma_{\text{ess}}(H_D)=0$.

In particular, $H_D$ possesses a compact resolvent if and only if 1-i holds.
\end{proposition}
To prove Proposition \ref{exab-out}, one uses \eqref{opposite}.
This essentially reduces the problem to the one considered in Proposition \ref{exabgen}.
\bigskip

We end this section with an example of the phenomenon mentioned in Remark 6.
On $R^+$ we give an example with $a_1=a_2=1$, and with $b_1$ and $b_2$ chosen appropriately so that   \eqref{intcondinfin} holds
for $a_1,b_1$ and $a_2,b_2$,  and such that
 $B_1(x)\equiv\int_0^xb_1(y)dy\ge B_2(x)\equiv\int_0^xb_2(y)dy$, but such that
\begin{equation}\label{b1}
\inf\sigma(-\frac12\frac{d^2}{dx^2}-b_1\frac d{dx})>0\ \ \text{and}\ \
\inf\sigma_{\text{ess}}(-\frac12\frac{d^2}{dx^2}-b_1\frac d{dx})=\infty,
\end{equation}
while
\begin{equation}\label{b2}
\inf\sigma(-\frac12\frac{d^2}{dx^2}-b_2\frac d{dx})=\inf\sigma_{\text{ess}}(-\frac12\frac{d^2}{dx^2}-b_2\frac d{dx})=0.
\end{equation}
 Let $b_1(x)=-x$ so that $B_1(x)=\int_0^xb_1(y)dy=-\frac{x^2}2$.
Then $\Omega^+(b_1,1)<\infty$ and $\hat \Omega^+(b_1,1)=0$, so \eqref{b1} holds.
It is not hard to construct a $b_2$ so that $B_2(x)<B_1(x)$, but such that for each positive integer $n$, there exists
an interval of length $n$ over which $b_2$ is identically 0. We will now show that $\Omega^+(b_2,1)
=\hat \Omega^+(b_2,1)=\infty$; thus, \eqref{b2} holds.
Using Theorem 1  and  the probabilistic representation in \eqref{exphitting}, we have for the diffusion corresponding to
$\frac12\frac{d^2}{dx^2}+b_2\frac d{dx}$ that
\begin{equation}\label{formula}
\frac1{8\Omega^+(b_2,1)}\le\sup\{\lambda\ge0: E_x\exp(\lambda\tau_0)<\infty\}\le\frac1{2\Omega^+(b_2,1)},\ x>0.
\end{equation}
Now for Brownian motion (that is, the driftless diffusion  corresponding to the operator $\frac12\frac{d^2}{dx^2}$)
on the interval $(0,n)$, one
has $E_x\exp(\lambda(\tau_0\wedge\tau_n))<\infty$, for $x\in(0,n)$, if and only if $\lambda$ is less than the first eigenvalue for
the operator  $-\frac12\frac{d^2}{dx^2}$ on $(0,n)$
with the Dirichlet boundary condition at 0 and $n$ \cite[chapter 3]{P}; that is, if and only if $\lambda<\frac{\pi^2}{2n}$.
Since the drift $b_2$ has intervals of length $n$ over which it vanishes, it follows by comparison with the Brownian motion that for the diffusion corresponding to
$\frac12\frac{d^2}{dx^2}+b_2\frac d{dx}$, if $x_n$ is chosen along such an interval,
then $E_{x_n}\exp(\lambda\tau_0)=\infty$, if $\lambda\ge\frac{\pi^2}{2n}$. Since the finiteness or infiniteness of the expectation
is independent of the starting point, it follows that in fact this holds for all $x>0$, not just for some $x_n$.
Since $n$ is arbitrary, it follows that
$\sup\{\lambda\ge0: E_x\exp(\lambda\tau_0)<\infty\}=0$. It then follows from \eqref{formula} that $\Omega^+(b_2,1)=\infty$, and then by the definition of
$\hat \Omega^+(b_2,1)$, also $\hat \Omega^+(b_2,1)=\infty$.

\section{Proofs of Theorems \ref{Th}-\ref{Th3} and of \eqref{exphitting}}
\bf\noindent Proof of Theorem \ref{Th}.\rm\
By Remark 2, it suffices to treat the case in which \eqref{intcondinfin} holds.
Extend $a$ and $b$ continuously from $[0,\infty)$ to $(-1,\infty)$. For each $l\in(-1,\infty)$, let
$H_D^{(l,\infty)}$ denote the corresponding self-adjoint diffusion operator on $(l,\infty)$ with the Dirichlet boundary condition
at $x=l$.
Consider the problem
\begin{equation}\label{posl}
\begin{aligned}
&\frac12(au')'+bu'+\lambda u=0,\  x\in(l,\infty);\\
&u>0,\ x\in(l,\infty).
\end{aligned}
\end{equation}
Let
\begin{equation}
\lambda_c(l)= \sup\{\lambda:\ \text{there is a solution to}\ \eqref{posl}\}.
\end{equation}
By the criticality theory of second-order elliptic operators,
there is a positive solution to the above equation for all $\lambda\le\lambda_c(l)$ \cite[chapter 4---section 3]{P} and
one has
$\inf\sigma(H_D^{(l,\infty)})=\lambda_c(l)$ \cite[chapter 4---section 10]{P}.
It follows from the  criticality theory that $\lambda_c(l)$ is right-continuous \cite[chapter 4---section 4]{P}.
However, in what follows we will need left continuity.
We claim that
\begin{equation}\label{cont}
\lambda_c(l) \ \text{is continuous in}\ l.
\end{equation}
We postpone the proof of \eqref{cont} until the end of the proof of Theorem \ref{Th}.
Note that any positive solution as above on $(l,\infty)$ with $l<0$ is also a positive solution
on $[0,\infty)$ and can be normalized by $u(0)=1$.
Note also that \eqref{pos} below always has a solution if $\lambda=0$.
From these facts, it follows that
if we consider the problem
\begin{equation}\label{pos}
\begin{aligned}
&\frac12(au')'+bu'+\lambda u=0,\ x\in(0,\infty);\\
&u>0,\ x\in(0,\infty);\\
&u(0)=1,
\end{aligned}
\end{equation}
then
\begin{equation}\label{char}
\inf\sigma(H_D)=\sup\{\lambda\ge0: \text{there is a solution to}\ \eqref{pos}\}.
\end{equation}
Thus, in order to prove \eqref{spec}, it suffices to prove the following proposition.

\begin{proposition}\label{ode}
Assume that \eqref{intcondinfin} holds.

\noindent i. For $\lambda>\frac1{2\Omega^+(b,a)}$, there is no solution to \eqref{pos};

\noindent ii. For $0<\lambda<\frac1{8\Omega^+(b,a)}$, there is a solution to \eqref{pos}.
\end{proposition}
The proof of part (i) of  Proposition \ref{ode}   is easy,  but the proof of part (ii) is
 nontrivial.
  The proof of the proposition  is given in  the next section.

Once \eqref{spec} is proved, one proves \eqref{essspec} as follows.
An old result of Persson \cite{Pe}, slightly modified to accommodate the case of a half-line, states that
\begin{equation}\label{essential}
\inf\sigma_{\text{ess}}(H_D)=\lim_{l\to\infty}(\inf\sigma(H_D^{(l,\infty)})).
\end{equation}
Letting
\begin{equation}\label{l}
\Omega^+_l(b,a)=\sup_{x>l}\left(\int_l^x\frac1{a(y)}\exp(-2B(y))dy\right)\left(\int_x^\infty\exp(2B(y))dy\right),
\end{equation}
 it follows by applying \eqref{spec} to $H^{(l,\infty)}_D$ that
 \begin{equation}\label{l-ineq}
 \frac1{8\Omega^+_l(b,a)}\le\inf\sigma(H_D^{(l,\infty)})\le\frac1{2\Omega^+_l(b,a)}.
\end{equation}
 We will show that
\begin{equation}\label{limitalpha}
\hat \Omega^+(b,a)=\lim_{l\to\infty}\Omega_l^+(b,a).
\end{equation}
Now \eqref{essspec} follows from \eqref{essential}, \eqref{l-ineq} and \eqref{limitalpha}.

We now prove \eqref{limitalpha}. From the definition of $\hat \Omega^+(b,a)$, one has for any $l>0$,
\begin{equation}
\begin{aligned}\label{onehand}
&\Omega^+_l(b,a)\ge\limsup_{x\to\infty}\left(\int_l^x\frac1{a(y)}\exp(-2B(y))dy\right)\left(\int_x^\infty\exp(2B(y))dy\right)\\
&=\limsup_{x\to\infty}\left(\int_0^x\frac1{a(y)}\exp(-2B(y))dy\right)\left(\int_x^\infty\exp(2B(y))dy\right)=\hat \Omega^+(b,a).
\end{aligned}
\end{equation}
On the other hand, for $n=1,2,\cdots$, there exist  $x_{0,n}$ and  $x_n$ with $x_{0,n}<x_n$ and
$\lim_{n\to\infty}x_n=\infty$, and such that
\begin{equation}\label{otherhand}
\begin{aligned}
&\limsup_{l\to\infty} \Omega^+_l(b,a)-\frac1n\le
\left(\int_{x_{0,n}}^{x_n}\frac1{a(y)}\exp(-2B(y))dy\right)\left(\int_{x_n}^\infty\exp(2B(y))dy\right)\\
&\le\left(\int_0^{x_n}\frac1{a(y)}\exp(-2B(y))dy\right)\left(\int_{x_n}^\infty\exp(2B(y))dy\right).
\end{aligned}
\end{equation}
Letting $n\to\infty$ in \eqref{otherhand} and again using the definition of $\hat \Omega^+(b,a)$, we obtain
 $\limsup_{l\to\infty}\Omega^+_l(b,a)\le\hat \Omega^+(b,a)$.
Now \eqref{limitalpha} follows from this and \eqref{onehand}.

We now return to prove \eqref{cont}.
As noted previously, we only need prove left-continuity.
Without loss of generality, we prove left-continuity at $l=0$.
From its definition, $\lambda_c$ is nondecreasing.
Let $\lambda_1<\lambda_2<\lambda_c(0)$.
It suffices to show  that for $\epsilon>0$ sufficiently small,
there is  a solution to \eqref{posl} with $l=-\epsilon$ and  some $\lambda\ge\lambda_1$.
By assumption, there is a solution to \eqref{posl} with $l=0$ and $\lambda=\lambda_2$.
Let $u$ be
 such a solution. Then $u(0^+)=\lim_{x\to0^+}u(x)$ and $u'(0^+)=\lim_{x\to 0^+}u'(x)$ exist and are finite.
 This is because any solution to \eqref{pos} must be a linear combination of
  $\Phi_1$ and $\Phi_2$, where $\Phi_1$ and $\Phi_2$ are two linearly independent solutions to
 $\frac12(au')'+bu'+\lambda u=0$. If $u(0^+)>0$, then solving the linear equation for $x<0$ using the boundary conditions
 $u(0^+)$ and $u'(0^+)$ at $x=0$, one can extend the solution $u$ a little bit to the left so that
it satisfies \eqref{posl} with $l=-\epsilon$ and $\lambda=\lambda_2$, completing the proof.

Assume now that $u(0^+)=0$. We will show that there exists a $\hat u$ which
 is a solution to \eqref{posl} with $l=0$ and $\lambda=\lambda_1$, and such that $\hat u(0)>0$.
Thus, from the previous argument, we can extend $\hat u$ a little bit to the left so that it
 satisfies \eqref{posl} with $l=-\epsilon$ and $\lambda=\lambda_1$, completing the proof.
Thus, it remains to show that such a $\hat u$ exists.
Let $\phi$ be a smooth compactly supported function on $R$ satisfying $\phi(0)=1$ and
$(\frac12(a\phi')'+b\phi'+\lambda_1\phi)(0)=0$.
Let $v=u+\delta \phi$, where $u$ is as above and  $\delta>0$. If $\delta$ is sufficiently small, then $v>0$ on $(0,\infty)$ and
$\frac12(av')'+bv'+\lambda_1v=-(\lambda_2-\lambda_1)u+\delta(\frac12(a\phi')'+b\phi'+\lambda_1\phi)<0$ on $(0,\infty)$.
Thus, $v$ is a sub-solution for \eqref{posl}
with $l=0$ and $\lambda=\lambda_1$, and $v(0)=\delta>0$.
We claim that there is a solution $\hat u$ to \eqref{posl} with $l=0$ and $\lambda=\lambda_1$, and with $\hat u(0)=\delta$.
Indeed, let $\hat u_n$ solve  $\frac12(a\hat u_n')'+b\hat u_n'+\lambda_1 \hat u_n=0$ in $(0,n)$, with $\hat u_n(0)=\delta$
and $\hat u_n(n)=0$. Then by the maximum principal, $\hat u_n$ is increasing in $n$ and
$\hat u_n\le v$; thus $\hat u\equiv\lim_{n\to\infty}\hat u_n$
exists and is the desired function.
\hfill$\square$

\medskip
\noindent \bf Proof of Theorem \ref{Th2}.\rm\
For $H_D$ on the entire line $R$, the result of Persson, given in \eqref{essential} for $R^+$, is
$$
\inf\sigma_{\text{ess}}(H_D)=\lim_{l\to\infty}\text{min}\left(\inf\sigma(H_D^{(l,\infty)}),\ \inf\sigma(H_D^{(-\infty,-l)})\right),
$$
where $H^{(-\infty,-l)}$ denotes the corresponding self-adjoint operator on $(-\infty,-l)$ with the Dirichlet boundary
condition at $x=-l$.
 Theorem \ref{Th2} follows from
this and the above proof of Theorem \ref{Th}.
\hfill $\square$

\medskip
\noindent\bf Proof of Theorem \ref{Th3}.\rm\
By the criticality theory of second-order elliptic operators \cite[chapter 4, sections 4 and 10]{P},
\begin{equation}\label{limitspec}
\inf\sigma(H_D)=\lim_{l\to\infty}(\inf\sigma(H_D^{(-l,\infty)}))=\lim_{l\to\infty}(\inf\sigma(H_D^{(-\infty,l)})).
\end{equation}
 Theorem \ref{Th} can be applied to $H_D^{(-l,\infty)}$. One simply lets $-l$ play the role played by
0 in Theorem \ref{Th}. (There is no need to change the definition of $B(x)=\int_0^x\frac ba(y)dy$, because
the lower limit 0 can be replaced by any $x_0$ without affecting the formulas.)
Thus, if $\int^\infty\frac1{a(x)}\exp(-2B(x))dx=\infty$, then defining
$$
\Omega^+_{-l}(b,a)=\sup_{x>-l}\left(\int_{-l}^x\frac1{a(y)}\exp(-2B(y))dy\right)\left(\int_x^\infty\exp(2B(y))dy\right),
$$
we have
\begin{equation}\label{spec-l}
\frac1{8\Omega^+_{-l}(b,a)}\le\inf\sigma(H_D^{(-l,\infty)})\le\frac1{2\Omega^+_{-l}(b,a)}.
\end{equation}
But
\begin{equation}\label{Omegalimit}
\lim_{l\to\infty}\Omega^+_{-l}(b,a)=\sup_{x\in R}
\left(\int_{-\infty}^x\frac1{a(y)}\exp(-2B(y))dy\right)\left(\int_x^\infty\exp(2B(y))dy\right).
\end{equation}
In the case that \eqref{recurrent} or \eqref{trans-} holds, Theorem \ref{Th3} follows from
\eqref{limitspec}-\eqref{Omegalimit}. The case that \eqref{trans+} holds is obtained from  the case
that \eqref{trans-} holds by interchanging the roles
of the positive and negative half-lines. For the  case that \eqref{trans+-} holds, one proceeds
as above in the case that \eqref{trans-} holds, but with $\Omega^+_{-l}(b,a)$ now defined
by
$$
\begin{aligned}
&\Omega^+_{-l}(b,a)=\\
&\sup_{x>-l}\left(\int_{-l}^xh^{-2}_{b,a}(y)\frac1{a(y)}\exp(-2B(y))dy\right)
\left(\int_x^\infty h^2_{b,a}(y)\exp(2B(y))dy\right).
\end{aligned}
$$
\hfill $\square$
\medskip

We end this section by proving that \eqref{exphitting} holds in the case that \eqref{intcondinfin} is in effect, that is, in the
case that $P_x(\tau_0<\infty)=1$. From \eqref{char}, it is enough to show that
\begin{equation}
\sup\{\lambda\ge0: \text{there is a solution to}\ \eqref{pos}\}=\sup\{\lambda\ge0: E_x\exp(\lambda\tau_0)<\infty\}.
\end{equation}
Assume first that $\lambda>0$ is such that there exists a solution to \eqref{pos} and let $u$ be a solution. Then
$u(X(t\wedge\tau_0))\exp(\lambda(t\wedge\tau_0))$ is a martingale \cite[chapter 2]{P}, and thus
\begin{equation}
E_xu(X(t\wedge\tau_0))\exp(\lambda (t\wedge\tau_0))=u(x).
\end{equation}
Letting $t\to\infty$, it follows from Fatou's lemma that
$E_x\exp(\lambda\tau_0)<\infty$.

Conversely, assume that $\lambda>0$ is such that $E_x\exp(\lambda\tau_0)<\infty$.
Let $\tau_n=\inf\{t\ge0:X(t)=n\}$, for $n>0$.
By the Feynman-Kac formula, $u_n(x)\equiv E_x(\exp(\lambda\tau_0);\tau_0<\tau_n)$ is the solution to the equation
\begin{equation}\label{pos-FK}
\begin{aligned}
&\frac12(au')'+bu'+\lambda u=0,\ x\in(0,n);\\
&u(0)=1,\
u(n)=0.
\end{aligned}
\end{equation}
By the maximum principal, $u_n$ is increasing in $n$, and \eqref{pos} will have a solution if and only if $\lim_{n\to\infty}u_n(x)<\infty$,
in which case $u_\infty(x)\equiv\lim_{n\to\infty}u_n(x)$ is the smallest solution to \eqref{pos}.
By the monotone convergence theorem and the assumption, we have $u_\infty(x)=E_x\exp(\lambda\tau_0)<\infty$.
Thus $\lambda$ is such that there is a solution to \eqref{pos}.

\section{Proof of Proposition \ref{ode}}
Let $\lambda>0$ and let $f_n$ be the unique solution to
\begin{equation*}
\begin{aligned}
&\frac12(af')'+bf'+\lambda f=0 \ \text{in}\ [0,n];\\
&f(0)=1, f(n)=0.
\end{aligned}
\end{equation*}
Integrating twice and using the boundary conditions gives
\begin{equation}\label{integraleq}
\begin{aligned}
&f_n(x)=1+c_n\int_0^x\frac1{a(y)}\exp(-2B(y))dy\\
&-2\lambda\int_0^xdy\frac1{a(y)}\exp(-2B(y))\int_0^ydzf_n(z)\exp(2B(z)),
\end{aligned}
\end{equation}
where
\begin{equation*}
c_n=\frac{-1+2\lambda\int_0^ndx\frac1{a(x)}\exp(-2B(x))\int_0^xdyf_n(y)\exp(2B(y))}{\int_0^n\frac1{a(x)}\exp(-2B(x))dx}.
\end{equation*}
Note that, by the maximum principle,  $f_n\ge0$ and  $f_n$ is nondecreasing in $n$.
Let $f_\infty\equiv\lim_{n\to\infty}f_n$.
Recall that by assumption,  $\int_0^\infty\frac1{a(x)}\exp(-2B(x))dx=\infty$.
Thus,
$c_\infty\equiv\lim_{n\to\infty}c_n=2\lambda\int_0^\infty f_\infty(x)\exp(2B(x))dx$.
Letting $n\to\infty$ in \eqref{integraleq} gives
\begin{equation}\label{f-infty}
f_\infty(x)=1+2\lambda\int_0^xdy\frac1{a(y)}\exp(-2B(y))\int_y^\infty dz f_\infty(z)\exp(2B(z)).
\end{equation}

By the maximum principle and the construction of $f_\infty$, either
$f_\infty$ is the smallest solution to \eqref{pos} or else $f_\infty=\infty$ and there are no  solutions to \eqref{pos}.
Using this characterization, we now proof the two parts of the proposition.

\noindent\it Proof of Part (i).\rm\
We will show that the solution $f_\infty$ of \eqref{f-infty} is equal to $\infty$
 if $\lambda>\frac1{2 \Omega(b,a)}$.
From \eqref{f-infty} it follows that $f_\infty$ is nondecreasing, and thus also that
\begin{equation}\label{contradiction}
f_\infty(x)\ge1+2\lambda\left(\int_0^x\frac1{a(y)}\exp(-2B(y))dy\right)\left(\int_x^\infty \exp(2B(y))dy\right)f_\infty(x).
\end{equation}
If there exists an $x$ for which $2\lambda(\int_0^x\frac1{a(y)}\exp(-2B(y))dy)(\int_x^\infty \exp(2B(y))dy)\ge1$, then
\eqref{contradiction} can not hold for such an $x$ unless $f_\infty(x)=\infty$.
Recalling the definition of $\Omega(b,a)$, we conclude that there is no finite solution to \eqref{f-infty}
if $\lambda>\frac1{2\Omega(b,a)}$.

\noindent \it Proof of Part (ii).\rm\ We will show that there is a finite solution to \eqref{f-infty} if $0<\lambda<\frac1{8\Omega(b,a)}$.
We assume that $\Omega(b,a)<\infty$ since otherwise there is nothing to prove. In particular then, we may assume
that $\int^\infty \exp(2B(z))dz<\infty$.
Fix $\lambda>0$ and define the operator
\begin{equation}
Tf(x)\equiv1+2\lambda\int_0^xdy\frac1{a(y)}\exp(-2B(y))\int_y^\infty dz f(z)\exp(2B(z)),
\end{equation}
operating on the domain $D_T\equiv\{f:f\ge0\ \text{and}\
\int^\infty f(z)\exp(2B(z))dz<\infty\}$. Note that, by assumption,
$1\in D_T$.
One can solve \eqref{f-infty} by iterations. Indeed,
it is clear that
 $T^n1$ is increasing in $n$ and that
$f_\infty=\lim_{n\to\infty}T^n1$, where $T^n$ denotes the $n$-th iterate of $T$.
Thus, to prove the existence of a finite solution to \eqref{f-infty} it is sufficient (and necessary) to show
that
\begin{equation}\label{suff}
\lim_{n\to\infty}T^n1<\infty.
\end{equation}
Define a norm by $||f||=\int_0^\infty f(x)\exp(2B(x))dx$.
We will prove \eqref{suff} by showing that
\begin{equation}\label{L1}
\lim_{n\to\infty}||T^n1||<\infty.
\end{equation}

Integrating by parts, we have
\begin{equation}
Tf=1+2\lambda S_1f+2\lambda S_2f,
\end{equation}
where
\begin{equation}
\begin{aligned}
&S_1f(x)=(\int_0^x\frac1{a(z)}\exp(-2B(z))dz)(\int_x^\infty f(z)\exp(2B(z))dz),\\
&S_2f(x)=\int_0^xdz f(z)\exp(2B(z))\int_0^zdt\frac1{a(t)}\exp(-2B(t)).
\end{aligned}
\end{equation}
Thus,
\begin{equation}\label{sumTn}
T^n1=1+\sum_{k=1}^n(2\lambda)^k(S_1+S_2)^k1.
\end{equation}
It is immediate from the definitions of $S_1$ and $\Omega(b,a)$ that $|S_11(x)|\le \Omega(b,a)$, and thus
\begin{equation}\label{1}
||S_11||\le \Omega(b,a)||1||.
\end{equation}
We will prove the following inequalities:
\begin{equation}\label{2}
||S_2f||\le \Omega(b,a)||f||;
\end{equation}
\begin{equation}\label{3}
||S^n_1S_2f||\le \Omega(b,a)(||S^{n-1}_1S_2f||+||S_1^nf||), \ n\ge1,
\end{equation}
where $S^0$ is defined to be the identity operator.
From \eqref{1}-\eqref{3}, it follows that
\begin{equation}\label{mix}
||S_{\delta_1}\cdots S_{\delta_k}1||\le(2\Omega(b,a))^k||1||,
\end{equation}
where $\delta_j=1$ or 2 for each $j=1, \cdots k$.
From \eqref{sumTn} and \eqref{mix} it follows that
\begin{equation}\label{Tn}
||T^n1||\le1+\sum_{k=1}^n (2\lambda)^k(2^k)(2\Omega(b,a))^k=1+\sum_{k=1}^n(8\lambda \Omega(b,a))^k.
\end{equation}
From \eqref{Tn} one concludes that \eqref{L1} holds if $\lambda<\frac1{8\Omega(b,a)}$.

We now prove \eqref{2} and \eqref{3}.
Integrating by parts, we have
\begin{equation}
\begin{aligned}
&||S_2f||=\int_0^\infty dx\exp(2B(x))\int_0^xdzf(z)\exp(2B(z))\int_0^zdt\frac1{a(t)}\exp(-2B(t))=\\
&-\left(\int_x^\infty\exp(2B(z))dz\right)\left(\int_0^xdzf(z)\exp(B(z))\int_0^zdt\frac1{a(t)}\exp(-2B(t))\right)\Big|_0^\infty\\
&+\int_0^\infty dx\left(\int_x^\infty  \exp(2B(z))dz\right)f(x)\exp(2B(x))\int_0^xdt\frac1{a(t)}\exp(-2B(t))\\
&\le \Omega(b,a)\int_0^\infty f(x)\exp(2B(x))dx=\Omega(b,a)||f||,
\end{aligned}
\end{equation}
proving \eqref{2}.

We now turn to \eqref{3}. We will write out the proof for $n=2$; the very same technique holds for general $n$.
We have
\begin{equation}\label{S112}
\begin{aligned}
&||S_1^2S_2f||=
\int_0^\infty dx\exp(2B(x))\left(\int_0^x\frac1{a(z)}\exp(-2B(z))dz\right)\times\\
&\left(\int_x^\infty dt \exp(2B(t))\int_0^tds\frac1{a(s)}\exp(-2B(s))
\int_t^\infty dl\exp(2B(l))S_2f(l)\right).
\end{aligned}
\end{equation}
Integrating by parts gives
\begin{equation}\label{S2}
\begin{aligned}
&\int_t^\infty dl\exp(2B(l))S_2f(l)=\\
&\int_t^\infty dl\exp(2B(l))\int_0^ldr f(r)\exp(2B(r))\int_1^rd\rho\frac1{a(\rho)}\exp(-2B(\rho))=\\
&-\left(\int_l^\infty \exp(2B(\nu))d\nu\right)
\left(\int_0^ldr f(r)\exp(2B(r))\int_0^rd\rho\frac1{a(\rho)}\exp(-2B(\rho))\right)\Big|_t^\infty+\\
&\int_t^\infty dl\left(\int_l^\infty d\nu \exp(2B(\nu))\right)f(l)\exp(2B(l))\int_0^ld\rho\frac1{a(\rho)}\exp(-2B(\rho))\le\\
&\left(\int_t^\infty \exp(2B(\nu))d\nu\right)
\left(\int_0^tdr f(r)\exp(2B(r))\int_0^rd\rho\frac1{a(\rho)}\exp(-2B(\rho))\right)+\\
&\int_t^\infty dl\left(\int_l^\infty d\nu \exp(2B(\nu))\right)f(l)\exp(2B(l))\int_0^ld\rho\frac1{a(\rho)}\exp(-2B(\rho)).
\end{aligned}
\end{equation}

Substituting \eqref{S2} in \eqref{S112} and using the definition of $\Omega(b,a)$, we obtain
\begin{equation*}
\begin{aligned}
&||S_1^2S_2f||\le\int_0^\infty dx\exp(2B(x))\left(\int_0^x\frac1{a(z)}\exp(-2B(z))dz\right)\times\\
&\Bigg(\int_x^\infty dt \exp(2B(t))\Big(\int_0^tds\frac1{a(s)}\exp(-2B(s))\Big)
\Big(\int_t^\infty d\nu\exp(2B(\nu))\Big)
\times\\
&\Big(\int_0^tdr f(r)\exp(2B(r))\int_0^rd\rho\frac1{a(\rho)}\exp(-2B(\rho))\Big)\Bigg)\\
&+\int_0^\infty dx\exp(2B(x))\Big(\int_0^x\frac1{a(z)}\exp(-2B(z))dz\Big)\times\\
&\Big(\int_x^\infty dt \exp(2B(t))
\big(\int_0^tds\frac1{a(s)}\exp(-2B(s))\big)\times\\
&\big(\int_t^\infty dl\big(\int_l^\infty d\nu\exp(2B(\nu))\big)f(l)\exp(2B(l))\int_0^lds\frac1{a(s)}\exp(-2B(s))\big)\Big)\\
&\le \Omega(b,a)\int_0^\infty dx\exp(2B(x))\left(\int_0^xdz\frac1{a(z)}\exp(-2B(z))\right)\times\\
&\Big(\int_x^\infty dt \exp(2B(t))
\big(\int_0^tdr f(r)\exp(2B(r))\int_0^rd\rho\frac1{a(\rho)}\exp(-2B(\rho))\big)\Big)\\
&+\Omega(b,a)\int_0^\infty dx\exp(2B(x))\Big(\int_0^x\frac1{a(z)}
\exp(-2B(z))dz\Big)\times\\
&\Big(\int_x^\infty dt \exp(2B(t))\big(\int_0^tds\frac1{a(s)}\exp(-2B(s))\big)
\big(\int_t^\infty dl f(l)\exp(2B(l))\big)\Big)\\
&=\Omega(b,a)||S_1S_2f||+\Omega(b,a)||S_1^2f||.
\end{aligned}
\end{equation*}
\hfill$\square$

\section{Application to Multi-Dimensional Diffusion Operators}
Consider the multi-dimensional diffusion operator
\begin{equation}\label{multi}
H_D=-\frac12\nabla\cdot a\nabla-a\nabla Q\cdot \nabla=-\frac12\exp(-2Q)\nabla \cdot a\exp(2Q)\nabla \ \text{ on}\ R^d,\ d\ge2,
\end{equation}
where $a=\{a_{i,j}\}_{i,j=1}^n\in C^1(R^d)$ is positive definite and $Q\in C^1(R^d)$.
One can realize $H_D$ as a non-negative, self-adjoint operator on $L^2(R^d,\exp(2Q)dx)$ via the
closure of the
 Friedrichs extension of the nonnegative
quadratic form
$$
Q_D(f,g)=\frac12\int_{R^d} \nabla fa\nabla g\exp(2Q)dx,
$$
 defined for $f,g\in C^1_0(R^d)$.
For $l>0$, let $B_l(0)\subset R^d$ denote the ball of radius $l$ centered at the origin, and let
 $H^l_D$ be the self adjoint operator on $R^d-\bar B_l(0)$ corresponding to $H_D$ with the Dirichlet boundary condition
 at $\partial B_l(0)$. More precisely, $H_D^l$ is the Friedrichs extension of the closure of the nonnegative quadratic form
$$
Q_D^l(f,g)=\frac12\int_{R^d-\bar B_l(0)} \nabla fa\nabla g\exp(2Q)dx,
$$
 defined for $f,g\in C^1_0(R^d-\bar B_l(0))$.
The result of Persson \cite{Pe} noted in section 3 gives
\begin{equation}\label{persson}
\inf\sigma_{\text{ess}}(H_D)=\lim_{l\to\infty}\inf\sigma(H_D^l).
\end{equation}
We will give upper and lower bounds on
 $\inf\sigma(H_D)$ and  $\inf\sigma(H_D^l)$ in terms of the corresponding infima for certain one-dimensional operators.
From \eqref{persson}, this will then also give upper and lower bounds on $\inf\sigma_{\text{ess}}(H_D)$.
Applying Theorem \ref{Th} to the one-dimensional operators will then yield explicit bounds on
$\inf\sigma(H_D)$ and $\inf\sigma_{\text{ess}}(H_D)$.

Letting $r=|x|$ and $\phi\in S^{d-1}$ denote spherical coordinates, let
\begin{equation}\label{har}
A_{\text{rad-har}}(r,\phi)=(\frac x{|x|}~a^{-1}(x)\frac x{|x|})^{-1}
\end{equation}
denote the representation in spherical coordinates of the reciprocal of the radially directed
quadratic expression $(\frac x{|x|}~a^{-1}(x)\frac x{|x|})$.
Let $Q_r(x)=\nabla Q(x)\cdot \frac x{|x|}$ denote the radial derivative of $Q$. Write $Q_r$ in spherical coordinates as $Q_r(r,\phi)$.
 For each $\phi\in S^{d-1}$,
define the one-dimensional diffusion operator $H_{\text{rad-har};\phi}$ on $R^+$ by
\begin{equation}\label{har}
\begin{aligned}
&H_{\text{rad-har};\phi}=\\
&-\frac12\frac d{dr}A_{\text{rad-har}}(r,\phi)\frac d{dr}-
A_{\text{rad-har}}(r,\phi)\frac{d-1}{2r}\frac d{dr}-
A_{\text{rad-har}}(r,\phi)Q_r(r,\phi)\frac d{dr}=\\
&-\frac12r^{1-d}\exp(-2Q(r,\phi))\frac d{dr}A_{\text{rad-har}}(r,\phi)~r^{d-1}\exp(2Q(r,\phi))\frac d{dr}
 \ \ \text{on}\ R^+.
\end{aligned}
\end{equation}

Let
$(\frac x{|x|}~a(x)\frac x{|x|})(r,\phi)$
denote the representation in spherical coordinates of the radially directed
quadratic expression $(\frac x{|x|}a(x)\frac x{|x|})$.
Let
\begin{equation}\label{avgdiff}
A_{\text{rad-avg}}(r)=\frac{\int_{S^{d-1}}(\frac x{|x|}~a(x)\frac x{|x|})(r,\phi)\exp(2Q(r,\phi))d\phi}
{\int_{S^{d-1}}\exp(2Q(r,\phi))d\phi}
\end{equation}
and
\begin{equation}\label{avgdrift}
Q_{r;\text{avg}}(r)=\frac{\int_{S^{d-1}}Q_r(r,\phi)\exp(2Q(r,\phi))d\phi}{\int_{S^{d-1}}\exp(2Q(r,\phi))d\phi}.
\end{equation}
Define the one-dimensional diffusion operator $H_{\text{rad-avg}}$ on $R^+$ by
\begin{equation}\label{avgop}
\begin{aligned}
&H_{\text{rad-avg}}=-\frac12\frac d{dr}A_{\text{rad-avg}}(r)\frac d{dr}-A_{\text{rad-avg}}(r)\frac{d-1}{2r}\frac d{dr}-
A_{\text{rad-avg}}(r)Q_{r;\text{avg}}(r)\frac d{dr}\\
&=-\frac12\exp(-2\beta(r))\frac d{dr}A_{\text{rad-avg}}(r)\exp(2\beta(r,\phi))~\frac d{dr}
 \ \ \text{on}\ R^+,\\
 &\text{where}\ \beta(r,\phi)=\frac{d-1}2\log r+\frac12\log\int_{S^{d-1}}\exp(2Q(r,\phi))d\phi.
 \end{aligned}
\end{equation}
Let $H_{\text{rad-har},\phi}^{(l,\infty)}$ and $H_{\text{rad-avg}}^{(l,\infty)}$ denote the corresponding operators
on $(l,\infty)$ with the Dirichlet boundary condition at $r=l$, as defined
in section 3.

\noindent \bf Remark 12.\rm\ Theorem \ref{Th} can be applied to the operators $H_{\text{rad-har},\phi}$ and $H_{\text{rad-avg}}$
even though their drifts are not continuous up to 0. Indeed, the theorem applies directly to
$H^{(l,\infty)}_{\text{rad-har},\phi}$ and $H^{(l,\infty)}_{\text{rad-avg}}$, for $l>0$,  and one has
$\inf\sigma(H_{\text{rad-har},\phi})=\lim_{l\to0}\inf\sigma(H^{(l,\infty)}_{\text{rad-har},\phi})$
and $\inf\sigma(H_{\text{rad-avg}})=$\newline$\lim_{l\to0}\inf\sigma(H^{(l,\infty)}_{\text{rad-avg}})$
\cite[chapter 4---sections 4 and 10]{P}. The one change that needs to be made is that $B$ should be
defined as $B(x)=\int_{x_0}^x\frac ba(y)dy$, for some $x_0>0$. (In the proof of Corollary \ref{Davies} below we use $x_0=1$.)

We will prove the following theorem.

\begin{theorem}\label{ineqavg}
$$
\begin{aligned}
&\inf_{\phi\in S^{d-1}}\inf\sigma(H_{\text{rad-har};\phi})
\le\inf\sigma(H_D)\le \inf\sigma(H_{\text{rad-avg}})\\
& \text{and}\\
&\inf_{\phi\in S^{d-1}}\inf\sigma(H^{(l,\infty)}_{\text{rad-har};\phi})\le
\inf\sigma(H^{(l,\infty)}_D)\le \inf\sigma(H^{(l,\infty)}_{\text{rad-avg}}).
\end{aligned}
$$
\end{theorem}

Applying Theorem \ref{Th} and \eqref{persson} to Theorem \ref{ineqavg}, the following corollary is immediate.
\begin{corollary}\label{corineqavg}
$$
\frac1{8\sup_{\phi\in S^{d-1}}\Omega^+(\beta_{\text{rad-har};\phi},A_{\text{rad-har}}(\cdot,\phi))}\le
\inf\sigma(H_D)\le\frac1{2 \Omega^+(\beta_{\text{rad-avg}},A_{\text{rad-avg}})};
$$
$$
\frac1{8\sup_{\phi\in S^{d-1}}\hat \Omega^+(\beta_{\text{rad-har};\phi},A_{\text{rad-har}}(\cdot,\phi))}\le
\inf\sigma_{\text{ess}}(H_D)\le\frac1{2\hat \Omega^+(\beta_{\text{rad-avg}},A_{\text{rad-avg}})},
$$
where
$$
\begin{aligned}
&\beta_{\text{rad-har};\phi}(r)=A_{\text{rad-har}}(r,\phi)\left(Q_r(r,\phi)+\frac{d-1}{2r}\right),\\
&\beta_{\text{rad-avg}}(r)=A_{\text{rad-avg}}(r)\left(Q_{r;\text{avg}}(r)+\frac{d-1}{2r}\right),
\end{aligned}
$$
and $\Omega^+$ and $\hat \Omega^+$ are as in Theorem \ref{Th}.
In particular, $\hat \Omega^+(\beta_{\text{rad-avg}},A_{\text{rad-avg}})=0$ is a necessary condition for $H_D$ to possess a compact resolvent
and \newline $\hat \Omega^+(\beta_{\text{rad-har};\phi},A_\text{rad-har}(\cdot,\phi))=0$, for all $\phi\in S^{d-1}$, is a sufficient condition.
\end{corollary}

\medskip

We give the following application of Corollary \ref{corineqavg}.
\begin{corollary}\label{Davies}
Let $H_D=-\frac12\nabla\cdot a\nabla$ on $R^d$, \ $d\ge1$.

\noindent i. If \ $\lim_{r\to\infty}\frac{\inf_{\phi\in S^{d-1}}A_{\text{rad-har}}(r,\phi)}{r^2}=\infty$, then
$\sigma_{\text{ess}}(H_D)=\emptyset$ and $H_D$ possesses a compact resolvent;

\noindent ii. If \  $\lim_{r\to\infty}\frac{A_{\text{rad-avg}}(r)}{r^2}=0$, then $\inf\sigma_{\text{ess}}(H_D)=0$;

\noindent iii. If   $\inf_{\phi\in S^{d-1}}A_{\text{rad-har}}(r,\phi)\ge\lambda r^2$, for large $r$, then
\newline $\inf\sigma_{\text{ess}}(H_D)\ge \frac{\lambda d^2}8$;

\noindent iv. If  $A_{\text{rad-avg}}(r)\le\Lambda r^2$, for large $r$, then
\newline $\inf\sigma_{\text{ess}}(H_D)\le \frac{\Lambda d^2}2$.
\end{corollary}
\bf\noindent Remark 13.\rm\
Let
$$
A_{\text{min}}(r)=\inf_{|v|=1, |x|=r}(va(x)v)\ \  \text{ and}\
\ A_{\text{max}}(r)=\sup_{|v|=1,|x|=r}(va(x)v),
 $$
  and note
 that
 $$
 \inf_{\phi\in S^{d-1}}A_{\text{rad-har}}(r,\phi)\ge A_{\text{min}}(r)\
 \  \text{and}\ \
 A_{\text{rad-avg}}(r)\le A_{\text{max}}(r).
$$
Parts (i)-(iii) of the above corollary, with
$\inf_{\phi\in S^{d-1}}A_{\text{rad-har}}(\cdot,\phi)$ replaced by  $A_{\text{min}}$ and
$A_{\text{rad-avg}}$ replaced by $A_{\text{max}}$
are originally due to Davies \cite{D}.
The use  of  $\inf_{\phi\in S^{d-1}}A_{\text{rad-har}}(\cdot,\phi)$
and  $A_{\text{rad-avg}}$
instead of $A_{\text{min}}$ and $A_\text{max}$ is a significant strengthening.
For instance, if for $|x|\ge1$, the radially directed vector $\frac x{|x|}$ is an eigenvector for $a(x)$ with
eigenvalue $\gamma(|x|)>1$, and all the other eigenvalues of $a(x)$ are equal to 1, then
for $r\ge1$, one has $A_{\text{rad-har}}(r,\phi)=\gamma(r)$ while  $A_{\text{min}}(r)=1$. A
two-dimensional example of such a diffusion matrix
is
$$
a(x)=\left(\begin{matrix} \gamma(|x|)\frac{x_1^2}{|x|^2}+\frac{x_2^2}{|x|^2}&
\frac{x_1x_2}{|x|^2}(\gamma(|x|)-1) \\ \frac{x_1x_2}{|x|^2}(\gamma(|x|)-1)&\gamma(|x|)\frac{x_2^2}{|x|^2}+
\frac{x_1^2}{|x|^2}  \end{matrix}\right).
$$
Switching the roles of the eigenvalues $\gamma(|x|)$ and 1 above, one has $A_{\text{rad-avg}}(r)=1$ while
$A_{\text{max}}(r)=\gamma(r)$.
A two-dimensional example of such a diffusion matrix
is
$$
a(x)=\left(\begin{matrix} \frac{x_1^2}{|x|^2}+\gamma(|x|)\frac{x_2^2}{|x|^2}&
\frac{x_1x_2}{|x|^2}(1-\gamma(|x|)) \\ \frac{x_1x_2}{|x|^2}(1-\gamma(|x|))&\frac{x_2^2}{|x|^2}+
\gamma(|x|)\frac{x_1^2}{|x|^2}  \end{matrix}\right).
$$
\medskip

\bf\noindent Proof of Corollary \ref{Davies}.\rm\
By the standard variational formula for $\inf\sigma(H^{(l,\infty)}_D)$, it follows that $\inf\sigma(H^{(l,\infty)}_D)$ is  nondecreasing in $a$,
and thus by \eqref{persson}, $\inf\sigma_{\text{ess}}(H_D)$ is also nondecreasing in $a$ .
Also from \eqref{persson}, it follows that $\inf\sigma_{\text{ess}}(H_D)$ does not depend on $\{a(x), 0<x\le l\}$, for any $l>0$.
In light of these facts, (i) and (ii) follow from (iii) and (iv).
Also, by the monotonicity in $a$, for the proof of (iii) we may assume that $A_{\text{rad-har}}(r,\phi)=\lambda r^2$, for large $r$,
 and  for the proof of (iv) we
may assume that $A_{\text{rad-avg}}(r)=\Lambda r^2$, for large $r$.

We consider (iii), the proof of (iv) following \it mutatis mutandi.\rm\
From Corollary \ref{corineqavg},
\begin{equation}\label{Daviesequ}
\inf\sigma_{\text{ess}}(H_D)\ge
\frac1{8\sup_{\phi\in S^{d-1}}\hat \Omega^+(A_{\text{rad-har}}(\cdot,\phi)\frac{d-1}{2r},A_\text{rad-har}(\cdot,\phi))}.
\end{equation}
For the pair of arguments of $\hat \Omega^+$ in \eqref{Daviesequ},
one has $\exp(2B(r))=\exp(\int_1^r\frac{d-1}sds)=r^{d-1}$.
Since we are assuming  that $A_{\text{rad-har}}(r,\phi)=\lambda r^2$ for large $r$,
we have $\int^\infty\frac1{A_{\text{rad-har}}(r,\phi)}r^{1-d}dr<\infty$; thus,
 \eqref{intcondfin} holds and
\begin{equation}\label{h-a}
h_{A_{\text{rad-har}}(\cdot,\phi)\frac{d-1}{2r},A_{\text{rad-har}}(\cdot,\phi)}(r)=
\int_r^\infty\frac1{A_{\text{rad-har}}(s,\phi)}s^{1-d}ds=\frac{r^{-d}}{\lambda d},\ \text{for large}\ r.
\end{equation}
Writing $h_\phi=h_{A_{\text{rad-har}}(\cdot,\phi)\frac{d-1}{2r},A_{\text{rad-har}}(\cdot,\phi)}$
to simplify notation,
for any $l>0$, one has from the definition of $\hat \Omega^+$,
\begin{equation}\label{Llatter}
\begin{aligned}
&\hat\Omega^+(A_{\text{rad-har}}(\cdot,\phi)\frac{d-1}{2r},A_\text{rad-har}(\cdot,\phi))=\\
&\limsup_{r\to\infty}\left(\int_l^rh_\phi^{-2}(s)
\frac1{A_{\text{rad-har}}(s,\phi)}s^{1-d}\right)\left(\int_r^\infty h_\phi^2(s)s^{d-1}ds\right)\\
&=\limsup_{r\to\infty}\left((h_\phi^{-1}(r)-h_\phi^{-1}(l)\right)
\int_r^\infty h_\phi^2(s) s^{d-1}ds.
\end{aligned}
\end{equation}
(In the original definition of $\hat \Omega^+$, $l$ above is replaced by 0, however using $l$ does not
 change the value of the expression.)
 Choosing $l$ sufficiently large and substituting
 $h_\phi(r)=\frac{r^{-d}}{\lambda d}$
 in  \eqref{Llatter}, one concludes  that\newline
$ \hat\Omega^+(A_{\text{rad-har}}(\cdot,\phi)\frac{d-1}{2r},A_\text{rad-har}(\cdot,\phi))=\frac1{\lambda d^2}$. Part (iii) now follows from
this and \eqref{Daviesequ}.
\hfill$\square$
\medskip

\noindent \bf Proof of Theorem \ref{ineqavg}.\rm\
We will prove the inequalities for $H_D$; the exact same method works for $H_D^{(l,\infty)}$.
The variational formula for $\inf\sigma(H_D)$  gives
\begin{equation}\label{RR}
\inf\sigma(H_D)=\inf\frac{\frac12\int_{R^d}\nabla fa\nabla f\exp(2Q)dx}
{\int_{R^d}f^2\exp(2Q)dx},
\end{equation}
where the infimum is over  $f\in C_0^1(R^d)$.

Using spherical coordinates $(r,\phi)$, and letting $\inf_{\text{radial}}$ denote the infimum
over radially symmetric functions $f\in C_0^1(R^d)$,
we have
from \eqref{RR},
\begin{equation}\label{estimate}
\begin{aligned}
&\inf\sigma(H_D)\le\inf_{\text{radial}}\frac{\frac12\int_{R^d}\nabla fa\nabla f\exp(2Q)dx}
{\int_{R^d}f^2\exp(2Q)dx}\\
&=\inf_{\text{radial}}\frac{\frac12\int_0^\infty (f'(r))^2\left(\int_{S^{d-1}}(\frac x{|x|}~a(x)\frac x{|x|})(r,\phi)\exp(2Q(r,\phi))d\phi\right)r^{d-1}dr}
{\int_0^\infty f^2(r)\left(\int_{S^{d-1}}\exp(2Q(r,\phi))d\phi\right)r^{d-1}dr}\\
&=\inf_{\text{radial}}\frac{\frac12\int_0^\infty (f'(r))^2A_{\text{rad-avg}}(r)\exp(2\beta(r))dr}
{\int_0^\infty(f^2(r)\exp(2\beta(r))dr},
\end{aligned}
\end{equation}
where
$\beta(r)=\frac{d-1}2\log r+\frac12\log\int_{S^{d-1}}\exp(2Q(r,\phi))d\phi$.
The infimum on the right hand side of \eqref{estimate} is the bottom of the spectrum of the
operator $H_{\text{rad-avg}}$  defined in \eqref{avgop}.
This gives the upper bound.

We now prove the lower bound.
By the Schwartz inequality,
\begin{equation}\label{schwartz}
\begin{aligned}
f^2_r(x)=(\nabla f(x)\cdot\frac x{|x|})^2\le (\nabla f(x)a(x)\nabla f(x))(\frac x{|x|}a^{-1}(x)\frac x{|x|}).
\end{aligned}
\end{equation}
Writing \eqref{schwartz} in polar coordinates and using the definition of $A_{\text{rad-har}}$, one has
\begin{equation}\label{schwartzapp}
(\nabla fa\nabla f)(r,\phi)\ge A_{\text{rad-har}}(r,\phi)f^2_r(r,\phi).
\end{equation}
By the variational formula, for any $g\in C^1_0(R^+)$ and $\phi\in S^{d-1}$,
\begin{equation}\label{Qphi}
\begin{aligned}
&\frac12\int_0^\infty A_{\text{rad-har}}(r,\phi)(g'(r))^2r^{d-1}\exp(2Q(r,\phi))dr\\
&\ge\inf\sigma(H_{\text{rad-har};\phi})\int_0^\infty g^2(r)r^{d-1}\exp(2Q(r,\phi))dr.
\end{aligned}
\end{equation}
From \eqref{schwartzapp} and \eqref{Qphi} one has for $f\in C^1_0(R^d)$,
\begin{equation}\label{lowerbound}
\begin{aligned}
&\frac12\int_{R^d}\nabla fa\nabla f\exp(2Q)dx
\ge\\
&\frac12\int_{S^{d-1}}\int_0^\infty A_{\text{rad-har}}(r,\phi)f_r^2(r,\phi)r^{d-1}\exp(2Q(r,\phi))drd\phi\\
&\ge\int_{S^{d-1}}\inf\sigma(H_{\text{rad-har};\phi})\int_{R^+}f^2(r,\phi)r^{d-1}\exp(2Q(r,\phi))drd\phi\\
&\ge\inf_{\phi\in S^{d-1}}\inf\sigma(H_{\text{rad-har};\phi})\int_{R^d}f^2\exp(2Q)dx.
\end{aligned}
\end{equation}
From \eqref{lowerbound} we conclude that
$$
\inf\sigma(H_D)=\inf\frac{\frac12\int_{R^d}\nabla fa\nabla f\exp(2Q(x))dx}
{\int_{R^d}f^2\exp(2Q)dx}\ge\inf_{\phi\in S^{d-1}}\inf\sigma(H_{\text{rad-har};\phi}).
$$
\hfill $\square$

\medskip
\bf\noindent Acknowledgement.\rm\ The author thanks Martin Kolb, a doctoral student
at the  University of  Kaiserslautern, for bringing to his attention the work of Muckenhoupt \cite{M} and for his
helpful comments.


\begin{thebibliography}{99}


\bibitem{D} Davies, E. B. $L\sp 1$ properties of second order elliptic operators, \emph{ Bull. London Math. Soc.}
\textbf{ 17} (1985),  417--436.

\bibitem{M} Muckenhoupt, B. Hardy's inequalities with weights, \emph{Studia Mathematica} \textbf{44} (1972), 31-38.


\bibitem{Pe}
Persson, A.,
Bounds for the discrete part of the spectrum of a semi-bounded Schr\"odinger operator,
\emph{Math. Scand.} \textbf{8}  (1960), 143-153.


\bibitem{P} Pinsky, R. G.,
\emph{Positive Harmonic Functions and Diffusion},
    Cambridge Studies in Advanced Mathematics
    \textbf{45}, Cambridge University Press,
      (1995).

\bibitem{P93} Pinsky, R. G., A new approach to the Martin boundary via diffusions conditioned to hit a compact set,
\emph{ Ann. Probab.} \textbf{21} (1993),  453--481.

\bibitem{RS2}
    Reed, M.  and Simon, B.,
    \emph{ Methods of Modern Mathematical Physics, II,  Fourier Analysis, Self Adjointness},
 Academic Press,
   New York, (1975).





\bibitem{RS4}
    Reed, M.  and Simon, B.,
    \emph{ Methods of Modern Mathematical Physics, IV,   Analysis of Operators},
 Academic Press,
   New York, (1975).



\end{thebibliography}
\end{document}